\documentclass[10pt]{amsart}

\usepackage{graphicx,amssymb,amsfonts,amsmath,amsthm,newlfont}
\usepackage{epsfig}
\usepackage{axodraw4j}
\usepackage{color,url}

\newtheorem{thm}{Theorem}
\newtheorem{conj}[thm]{Conjecture}
\newtheorem{ques}[thm]{Question}
\newtheorem{lem}[thm]{Lemma}
\newtheorem{prop}[thm]{Proposition}
\newtheorem{propdef}[thm]{Proposition-Definition}
\newtheorem{remark}[thm]{Remark}
\newtheorem{cor}[thm]{Corollary}

\theoremstyle{definition}
\newtheorem{defn}[thm]{Definition}
\newtheorem{example}[thm]{Example}

\newcommand{\A}{\mathbb{A}}

\newcommand{\Lef}{\mathbb{L}}
\newcommand{\M}{\mathbb{M}}
\newcommand{\Q}{\mathbb Q}
\newcommand{\Z}{\mathbb Z}
\newcommand{\R}{\mathbb R}
\newcommand{\F}{\mathbb F}
\newcommand{\Pro}{\mathbb P}

\newcommand{\gr}{\mathrm{gr}}

\newcommand{\q}{/\!\!/}

\title[A K3 in $\phi^4$]{A  K3 in $\phi^4$}

\author{Francis Brown and Oliver Schnetz}

\begin{document}
\begin{abstract}
 Inspired by Feynman integral computations in quantum field theory, Kontsevich conjectured in 1997 that the number of points of graph hypersurfaces over a finite
 field $\F_q$ is a (quasi-) polynomial in $q$.  Stembridge verified this  for all graphs with  $\leq12$ edges, but  in 2003 Belkale and Brosnan showed  that the counting
 functions are of general type for  large graphs. In this paper we 
 give  a sufficient combinatorial criterion  for a graph to have   polynomial  point-counts,   and construct some explicit counter-examples to Kontsevich's conjecture
 which are in $\phi^4$ theory. Their   counting functions  are given modulo $pq^2$   ($q=p^n$) by a modular form arising from a certain singular K3 surface. 
\end{abstract}
\maketitle

\section{Introduction}
We first recall the definition of graph hypersurfaces and the  history of the point-counting problem, before explaining its relevance to Feynman integral calculations
in perturbative Quantum Field Theory.

\subsection{Points on graph hypersurfaces}
Let $G$ be a connected graph, possibly with multiple edges and self-loops (an edge whose endpoints coincide). 
The graph polynomial of $G$  is defined by associating  a variable $\alpha_e$, known as a Schwinger parameter, to every edge $e$ of $G$ and  setting
\begin{equation} \label{intropsiG}
\Psi_G(\alpha) = \sum_{T\subset G} \prod_{e\notin T} \alpha_e \in \Z[\alpha_e]\ , 
\end{equation}
where the sum is over all spanning trees $T$ of $G$ (connected subgraphs  meeting every vertex of $G$ which have no loops). These polynomials go back to the work of  Kirchhoff
in relation to the study of 
currents in electrical circuits \cite{Kir}.

The projective graph hypersurface $X_G$ is defined to be the zero locus of $\Psi_G$ in projective space $\Pro^{N_G-1}$, where $N_G$ is the number of edges of $G$
(although from \S\ref{sectresults} onwards,   $X_G$ will denote the zero locus in affine space $\A^{N_G}$).
It is highly singular in general. For any prime power $q$, let $\F_q$ denote the field with $q$ elements, and consider the point-counting function:
$$[X_G]_q : q \mapsto \# X_G(\F_q) \in \mathbb{N}\cup\{0\}.$$
In 1997, Kontsevich informally conjectured that this function might be  polynomial in $q$ for all graphs. 
This question was studied by Stanley, Stembridge   and  others, and in particular was proved for all graphs with at most twelve edges \cite{Stem}.
A dual statement was proved for various families of graphs obtained by deleting trees in complete graphs \cite{Sta}, \cite{CY}. But in \cite{BB}, contrary to expectations,
Belkale and Brosnan used Mn\"ev's universality theorem to  prove    that the $[X_G]_q$ are of general type  in the following precise sense.

\begin{thm} (Belkale-Brosnan). For every scheme $Y$ of finite type  over $\mathrm{Spec}\,\,  \Z$, there exist finitely many polynomials $p_i \in \Z[q]$  and graphs $G_i$ such that
$$ s [Y]_q = \sum_{i} p_i [X_{G_i}]_q\ ,$$
where $[Y]_q$ denotes the point-counting function on $Y$, and $s$ is a product  of terms of the form $q^n-q$,  where $n>1$. In particular,
$[X_G]_q$ is not always polynomial.
\end{thm}

This does not imply that the point-counting functions 
$[X_{G_i}]_q$ themselves are arbitrary.  The methods of \cite{BlochJ} \S4, for example, imply strong constraints on
$[X_{G}]_q$.
Moreover, Belkale and Brosnan's method constructs graphs $G_i$ with very large numbers of edges (\cite{BB}, Remark 9.4), and no explicit counter-example was known until recently,
when Doryn \cite{DorynPoints}
and Schnetz \cite{SchnetzFq} independently constructed graphs which are quasi-polynomial (i.e., which become polynomial only  after a finite extension of the  base field
or the exclusion of exceptional primes). 
It has since been  hoped that various  `physicality' constraints  on $G$ might be sufficient to ensure the validity of Kontsevich's conjecture in this  slightly weaker sense. 
However, the modular counter-examples we construct below show that this hope is completely false.

\subsection{Feynman integrals and motives}  The point-counting problem has its origin in the question of determining the arithmetic content of perturbative quantum field theories.
For this, some convergency conditions are required on the graphs. A connected graph $G$ is said to be primitive-divergent if:
\begin{eqnarray} \label{defnprimdiv}  N_G  & = & 2h_G  \\
 N_{\gamma} &> &2 h_{\gamma} \quad \hbox{ for all strict subgraphs } \gamma \subsetneq G \ ,\nonumber
 \end{eqnarray}  where $h_{\gamma}$ denotes the number of loops (first Betti number) and 
$N_{\gamma}$ the number of edges in a graph 
$\gamma$.
In this case, the residue of $G$  is defined by the absolutely convergent projective integral \cite{Wein}, \cite{B-E-K}
\begin{equation} \label{IG}
I_G = \int_{\sigma} \frac{\Omega_{N_G}}{\Psi_G^2}\ ,
\end{equation}
where $\sigma=\{(\alpha_1:\ldots :\alpha_{N_G}) \subset \Pro^{{N_G}-1}(\R) : \alpha_i \geq 0 \}$ is the real coordinate simplex in projective space, and
$\Omega_{N_G}=\sum_{i=1}^{N_G} (-1)^i d\alpha_1 \ldots \widehat{d\alpha_i} \ldots d\alpha_{N_G}$. This defines a map from the set of primitive-divergent graphs to positive real numbers.
It is important to note that the quantities $I_G$  are renormalization-scheme independent.  
We say that $G$ is in $\phi^4$ theory if every vertex of $G$ has degree at most four.
Even in this case, the numbers $I_G$ are very hard to evaluate, and known analytically for only a handful of graphs.
Despite the difficulties in computation, the remarkable fact was observed by Broadhurst, Kreimer  \cite{BK}, and later  Schnetz \cite{SchnetzCensus},
that every graph whose period is computable (either analytically or
numerically to high precision) is consistent with being a multiple zeta value. This was the original motivation for Kontsevich's question.
\vskip2ex

The algebraic  approach to this problem comes from  the observation  that the numbers $I_G$ are periods in the sense of algebraic geometry.
To make this precise,   the integrand of $(\ref{IG})$ defines a cohomology class in $H^{N_G-1} (\Pro^{N_G-1} \backslash X_G)$, and the domain of integration
a relative homology class in $H_{N_G-1} (\Pro^{N_G-1}, B)$ where $B=V(\prod_{i=1}^{N_G} \alpha_i)$, which contains the boundary of the simplex $\sigma$.
Thus   as a first approximation, one could consider the relative mixed Hodge structure
\begin{equation}\label{HN}  H^{{N_G}-1}(\Pro^{N_G-1}\backslash X_G, B\backslash (B\cap X_G))\ .
\end{equation}
For technical reasons related to the fact that $\sigma$ meets $X_G$ non-trivially, the integral $I_G$ is not in fact a period of $(\ref{HN})$. One of the main constructions of
\cite{B-E-K}  is to blow up boundary components of $B$ to obtain a slightly different relative mixed Hodge structure called  the graph motive $M_G$.
The integral $I_G$ is now a period of $M_G$.
If $M_G$ is of mixed Tate type (its weight graded pieces are of type $(p,p)$) and satisfies some ramification conditions, then by standard conjectures on mixed Tate motives (now  proved  \cite{BrAn}),
it should follow that  the period $I_G$ is  a multiple zeta value. 

Although not explicitly stated in \cite{B-E-K}, it follows from the geometry underlying their construction and the relative cohomology  spectral sequence that $M_G$ 
is controlled by the absolute mixed Hodge structures 
$H^i(\Pro^i \backslash X_{\gamma})$, where $\gamma$ ranges over all minors (subquotients) of $G$. Thus the simplest way in which the period $I_G$ could be a multiple zeta value is 
if  the  mixed Hodge structure $M_{\gamma}$ were entirely of Tate type, or, even stronger, if  $H^{\bullet} (\Pro^i \backslash X_{\gamma})$ were of Tate type
in all cohomological dimensions, for all minors $\gamma$ of $G$. 
To simplify matters  further, one  can  ask  the somewhat  easier  question of  whether the  Euler characteristics  of the $X_{\gamma}$'s are of Tate type.
In this way, one is  led to consider the class of $X_G$ in the Grothendieck ring of varieties $K_0(\hbox{Var}_k)$ and ask if it is a polynomial in the Lefschetz motive $\Lef=[\A^1_k]$.
This is surely the reasoning behind Kontsevich's original question, although it was formulated almost ten years before $M_G$ was defined.
Note, however,  that there is \emph{a priori} no way to construe information about $I_G$ from  the Grothendieck class $[X_G]$.

\subsection{Results and contents of the paper} \label{sectresults}
We begin in $\S2$ by reviewing some algebraic and combinatorial properties of graph polynomials.  In $\S3$, we discuss implications for the class of the affine
graph hypersurface $[X_G]$ in the Grothendieck ring of varieties $K_0(\hbox{Var}_k)$, where $k$ is a field. The first observation is the following:

\begin{prop}  Let $G$ be any graph satisfying $h_G\leq N_G-2$. Then  there is an invariant  $c_2(G) \in K_0(\hbox{Var}_k)/ \Lef$ 
such that
\begin{equation}\label{introc2def}
[X_G] \equiv c_2(G) \Lef^2 \mod \Lef^3\ .
\end{equation}
If $G$ has a three-valent vertex, then   $c_2(G)$ has a simple representative   in $K_0(\hbox{Var}_k)$ as the class of the intersection of two explicit affine hypersurfaces. 
\end{prop}

For primitive-divergent graphs  $(\ref{defnprimdiv})$ this intersection satisfies a Calabi-Yau property in the sense that, after projectifying, the total
degree is exactly one greater than the dimension of the ambient projective space.

The  class $c_2(G)$ is much more tractable than the full class $[X_G]$.  In order to exploit its combinatorial properties, we require  the Chevalley-Warning theorem \ref{thmCW}
on the point-counts   of polynomials of small degree modulo $q$. 
We hope that this theorem  lifts to  the Grothendieck ring  under some conditions on $k$  (\S3.3), but 
since this is   unavailable, we are  forced to  pass to point-counts modulo $q$. Thus, denoting the corresponding counting functions by $[.]_q$, equation $(\ref{introc2def})$ gives
$$[X_G]_q \equiv c_2(G)_q q^2 \mod q^3\ ,$$
where  $c_2(G)_q$ is a map from prime powers $q$ to 
$\Z/q\Z$.

In $\S3$ we explain how to compute the invariant $c_2(G)$ by a simple  algorithm (`denominator reduction')  which 
reduces the problem to counting points on hypersurfaces of  smaller and smaller  dimension. This  is the key to constructing non-Tate counter-examples,
and stems from the following observation:
\begin{thm} Let $G$ be primitive divergent with at least five edges $e_1,\ldots,e_5$. Then
$$c_2(G)_q \equiv -[{}^5\Psi_{G}(e_1,\ldots, e_5)]_q \mod q$$
where  ${}^5\Psi_G(e_1,\ldots, e_5)$ is  the 5-invariant \cite{BrFeyn} of those edges.  
\end{thm}

The 5-invariant is a certain resultant of polynomials derived from the graph polynomials of minors of $G$, and it follows from this theorem that 
the invariant $c_2(G)_q$ can be computed inductively by taking iterated resultants. This uses  the Chevalley-Warning theorem in an essential way  to kill parasite terms. 

In $\S4$ we use this result to deduce the following properties of $c_2(G)_q$ in the case when $G$ is 
primitive-divergent graph in $\phi^4$ theory:
 \begin{enumerate}
\item If $G$ is two-vertex reducible then $c_2(G)_q \equiv 0 \mod q$.
 \item If $G$ has weight-drop (in the sense of \cite{WD}), then $c_2(G)_q \equiv 0 \mod q$.
  \item If $G$ has vertex-width $\leq 3$, then $c_2(G)_q \equiv c \mod q$ for some $c\in \Z$.
 \item $c_2(G)_q$ is invariant under double triangle reduction.
 \end{enumerate} 
 For the definitions of these terms, see \S4. All of these properties have some bearing on the residue $I_G$ (\cite{BrFeyn}, \cite{WD}). A further property is conjectural:
\begin{conj} $c_2(G)$ is invariant under the completion relation (\cite{SchnetzCensus}, \cite{SchnetzFq}).
\end{conj} 

 This is implied by  the following stronger conjecture (see \cite{SchnetzFq}, remark 2.10 (2)):
\begin{conj} If $I_{G_1}=I_{G_2}$ for two graphs $G_1$, $G_2$ then $c_2(G_1)=c_2(G_2)$.
\end{conj} 

In short, the invariant $c_2(G)$ detects all the known  qualitative features of the residue $I_G$, but is much easier to compute. 
Intuitively, $c_2(G)_q$ should be closely related to the action of Frobenius on the  framing of $M_G$, i.e., the smallest subquotient motive of $M_G$ which is spanned by the
Feynman differential form $\Omega_{N_G}\Psi_G^{-2}$.

In $\S5$ we review the notion of vertex-width, which is a measure of the  local connectivity of a graph, and prove Kontsevich's conjecture for an infinite family of graphs.
Note that this result is valid in the Grothendieck ring $K_0(\hbox{Var}_k)$.
\begin{thm} \label{thmvw3intro} Let $G$ have vertex-width at most 3. Then $[\Psi_G]$ is a polynomial in $\Lef$.
\end{thm}
This family of graphs contains  almost all the physically interesting cases at low loop orders.
It was proved in \cite{BrFeyn} that a variant of the motive $M_G$ is mixed Tate in this case, but the proof we give here is elementary and gives an effective way to compute the 
polynomial $[\Psi_G]$ by induction over the minors of $G$. It also enables one to compute the Grothendieck classes of any infinite family of graphs obtained by inserting triangles
into a known graph (a problem raised in \S13.2 of \cite{B-E-K}).
In \S\ref{sectwheels} and  \S\ref{sectZigs}, we carry this out for the 
wheels  (also computed independently in   \cite{DorynPoints}) and zig-zag graphs.  
These  are the only two families  of graphs  for which a 
formula for the residue $I_G$ is known, or conjectured.

The motivation for such computations is the hope that they will  give combinatorial insight into the full structure of the motive $M_G$, and ultimately the action of the motivic
Galois group, which would yield a lot of information about the periods. Currently there is not a single example where this has been successfully carried out.

 \begin{figure}[ht!]
 \begin{center}
    \leavevmode
    \epsfxsize=7.0cm \epsfbox{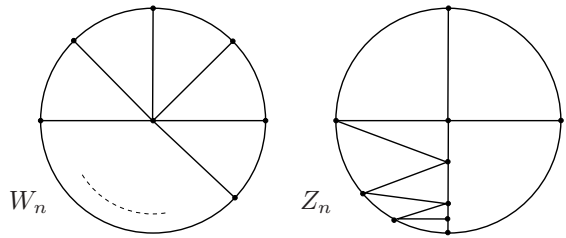}
  \put(-210,10){{$W_n$}}  \put(-100,10){{$Z_n$}}
  \end{center}
  \caption{The wheels with spokes  (left), and zig-zags (right).}
\end{figure}

In $\S6$ we construct explicit counter-examples to Kontsevich's conjecture by computing the  $c_2$-invariants of some graphs $G$ of vertex width 4.
The denominator reduction algorithm  reduces  $c_2(G)_q$ down to a determinant of graph polynomials of small degree derived from $G$.
By a series of manipulations  one can extract    a polynomial which defines a surface of degree 4 in $\Pro^3$, whose minimal desingularization $X$  is a K3 surface.
In $\S7$, we  show that this surface has N\'eron-Severi group  of maximal rank, and that its Picard lattice has  discriminant $-7$. This proves that $X$ is a singular
$K3$ surface, which have been classified by Shioda and Inose \cite{SI}. The modularity of such surfaces is known by \cite{Liv}, and in this case $H^2_{tr}(X)$ is a submotive of
the symmetric square of the first cohomology group of the elliptic curve:
$$E_{49A1}: \qquad y^2+xy = x^3-x^2-2x-1\ ,$$
which has complex-multiplication by $\Q(\sqrt{-7})$.  We write down the modular form of 
weight 2 and level 49 whose coefficients give the point-counts on $E_{49A1}$. Its symmetric square is given by the following product of Dedekind $\eta$-functions:
\begin{equation} \big(\eta(z) \eta(z^7) \big)^3 \ , \label{introeta}
\end{equation} 
which defines a cusp form of weight 3 and level 7.
\begin{thm} \label{introthmcounterex}  Let $q=p^n$. There exists a non-planar  primitive-divergent graph  in $\phi^4$ 
theory with 8 loops, 16 edges, and vertex-width 4 such that 
$$c_2(G)_q \equiv - a_q^2 \equiv -b_q \mod p$$
where  $q+1- a_q = [E_{49A1}]_q  $ is the number of points on $E_{49A1}(\F_q)$, and $b_q$ is  the coefficient of $z^q$ in $(\ref{introeta}).$  In particular, 
\begin{equation}\label{1}
[X_G]_q \equiv -a^2_q q^2 \equiv -b_q q^2 \mod pq^2
\end{equation}
cannot be a polynomial in $q$. Furthermore,
there exists a  {\bf planar} primitive-divergent graph in $\phi^4$ theory with 9 loops, 18 edges, and vertex-width 4 with the same property.
\end{thm}
The fact that our counter-examples have vertex-width 4 shows  that theorem \ref{thmvw3intro} cannot be improved.
In fact, within $\phi^4$ theory the 8-loop graph  is the smallest modular  counter-example to Kontsevich's conjecture 
\cite{SchnetzFq}.
\vskip2ex

Both authors wish to thank D. Broadhurst, D. Kreimer, H. Esnault,  M. Sch\"utt,  K. Yeats and especially S. Bloch, and thank Humboldt University, Berlin,   for hospitality.
Francis Brown is supported by ERC grant 257638.

\section{Graph polynomials}
Throughout this paper, $G$ will denote  any graph, possibly with multiple edges and self-loops. A subgraph of $G$ will be a subset of edges of $G$.
In this section, we make no assumptions about the primitive-divergence or otherwise of $G$.

\subsection{Matrix representation}  We recall some basic results from \cite{BrFeyn}. We will use the following matrix representation for the   graph polynomial. 

\begin{defn} Choose an orientation on the edges of $G$, and for every edge $e$ and vertex $v$ of $G$, define the incidence matrix:
$$(\mathcal{E}_G)_{e,v} = \left\{
                           \begin{array}{rl}
                             1, & \hbox{if the edge } e \hbox{ begins at } v \hbox{ and does not end at } v,\\
                             -1, & \hbox{if the edge } e  \hbox{ ends at } v \hbox{ and does not begin at } v,\\
                             0, & \hbox{otherwise}.
                           \end{array}
                         \right.
 $$
Let $A$ be the diagonal matrix with entries  $\alpha_e$, for $e \in E(G)$, and set
$$\widetilde{M}_G=\left(
  \begin{array}{c|c}
    A  & \mathcal{E}_G  \\
    \hline
  {-}\mathcal{E}_G^T&  0  \\
  \end{array}
\right)
$$
where the first $N_G$ rows and columns are indexed by the set of edges of $G$, and the remaining $v_G$ rows and columns are indexed by the set of vertices of $G$, in some order.
The matrix $\widetilde{M}_G$ has corank  $\geq 1$.  Choose any vertex of $G$ and let $M_G$ denote the square $(N_G+v_G-1)\times (N_G+v_G-1)$ matrix obtained from it
by deleting the row and column indexed by this vertex.
\end{defn}

It follows from the matrix-tree theorem that the graph polynomial satisfies (\S2.2 in \cite{BrFeyn})
$$\Psi_G=\det (M_G)\ .$$
If $G$ has at  least two  components $G_1$ and $G_2$,  then by permuting rows and columns $\widetilde{M}_G$ can be written as a direct sum of $\widetilde{M}_{G_1}$ and
$\widetilde{M}_{G_2}$. In this case the graph polynomial   $\Psi_G$ vanishes, since
$\mathrm{corank}\, \widetilde{M}_G = \mathrm{corank} \, \widetilde{M}_{G_1} +\mathrm{corank}\, \widetilde{M}_{G_2}\geq 2.$

\begin{defn}  Let $I,J, K$ be subsets of the set of edges of $G$ which satisfy $|I|=|J|$. Let  
$M_G(I,J)_K$ denote the matrix obtained from $M_G$ by removing the rows (resp. columns) indexed by the set $I$ (resp.\ $J$) and setting $\alpha_e=0$ for all $e\in K$. 
Let 
\begin{equation} \label{Dogsondefn} \Psi_{G,K}^{I,J}=\det M_G(I,J)_K\ . \end{equation}
We shall perpetuate the anachronism of \cite{BrFeyn} by referring  to these polynomials as Dodgson polynomials.
\end{defn}

It is clear that $\Psi_{G,\emptyset}^{\emptyset,\emptyset} = \Psi_G$, and $\Psi^{I,J}_{G,K} = \Psi^{J,I}_{G,K}$ because $I,J\subset E(G)$.
If $K=\emptyset$, we will often drop it from the notation. We also write $\Psi^I_{G,K}$ as a shorthand for $\Psi^{I,I}_{G,K}$.

Since the matrix $M_G$ depends on various choices, the polynomials $\Psi^{I,J}_{G,K}$ are only well-defined up to sign. In what follows, for any graph $G$, we shall fix a particular
matrix $M_G$ and this will fix all the signs in the polynomials $\Psi^{I,J}_{G,K}$ too.

\begin{prop} \label{propspanningtrees} The monomials which occur in  $\Psi^{I,J}_{G,K}$ have coefficient $\pm 1$, and are precisely the monomials which occur in both
$\Psi^{I,I}_{G,J\cup K}$ and $\Psi^{J,J}_{G,I \cup K}$.
\end{prop}
\begin{proof} See \cite{BrFeyn}, proposition 23, $\S2.3$.
\end{proof} 

\begin{defn}
If $f=f_1+f^1 \alpha_1$ and $g=g_1+g^1 \alpha_1$ are polynomials of degree one in $\alpha_1$, recall that their resultant is defined by:
\begin{equation} \label{resultantdef}
[f,g]_{\alpha_1}= f^1g_1-f_1g^1\ .
\end{equation}
\end{defn} 
We now state some identites between Dodgson polynomials which will be used in the sequel. The proofs can be found in (\cite{BrFeyn}, \S2.4-2.6).

\subsection{General identities}
The first set of identities only use
 symmetries of the matrix $M_G$, and therefore hold for any graph $G$.
\begin{enumerate}
\item \emph{The 
contraction-deletion formula}.  It is clear from its definition that $\Psi^{I,J}_{G,K}$ is linear in every Schwinger variable $\alpha_e$. When the index of $e$ is larger than
all elements of $I$ and $J$ (in general there are signs), we  can   write:
$$\Psi^{I,J}_{G,K} = \Psi^{I e, Je}_{G,K} \alpha_e + \Psi^{I,J}_{G,Ke}\ .$$
The  contraction-deletion relations state that
  $$\Psi^{Ie,Je}_{G,K}=\Psi^{I,J}_{G\backslash e, K} \, \hbox{ and }  \Psi^{I,J}_{G,Ke} = \Psi^{I,J}_{G\q e, K}\ , $$
   where $G\backslash e$ is the  graph obtained by deleting the edge $e$ (but not its endpoints), and $G\q e$ denotes the graph obtained
 by contracting the edge $e$ (and identifying its two endpoints).
Note that we define the contraction of a self-loop to be  the zero graph $0$, for which we set  $\Psi_0=0$.
 \vspace{0.05in}
\item \emph{Dodgson-type identities}. Let $I,J$ be two subsets of edges of $G$ such that $|I|=|J|$ and let $a,b,x\notin I \cup J
\cup K$. Then the first (Dodgson) identity is:
$$ \big[\Psi^{I,J}_{G,K}, \Psi^{Ia,Jb}_{G,K} \big]_x = \Psi^{Ix,Jb}_{G,K} \Psi^{Ia,Jx}_{G,K}\ .$$
 Let $I,J$ be two subsets of edges of $G$ such that $|J|=|I|+1$ and let $a,b,x\notin I\cup J
\cup K$. Then the second identity is:
$$ \big[\Psi^{Ia,J}_{G,K}, \Psi^{
Ib,J}_{G,K} \big]_x = \pm \Psi^{Ix,J}_{G,K} \Psi^{Iab,Jx}_{G,K}\ . $$
\vspace{0.05in}
\end{enumerate}
Note that $\Psi^{I,I}_{G,K}=\Psi_{G\backslash I \q K}$. A graph obtained by contracting and deleting edges of $G$ will be called a minor of $G$.
\subsection{Graph-specific identities} \label{sectGraphSpecIds}
The second set of identities depend on the particular combinatorics of a graph $G$, and follow from proposition 
$\ref{propspanningtrees}$ together with the fact that
$\Psi_{G\backslash I}=0$ if $I$ contains the set of  all edges which meet  a given vertex
(because $G\backslash I$ is not connected), and $\Psi_{G\q K} =0$ if $h_1(K)>0$ 
(because contracting edges of $K$ one by one leads to  the contraction of a self-loop).
\begin{enumerate}

\item \emph{Vanishing property for vertices}. Suppose that $E=\{e_1,\ldots, e_k\}$ is the set of  edges which are adjacent to a given  vertex of $G$. Then
$$\Psi^{I,J}_{G,K} =0 \qquad  \hbox{ if } \qquad  E\subset I \quad \hbox{ or } \quad E \subset J\ .$$

\item \emph{Vanishing property for loops}. Suppose that $E=\{e_1,\ldots, e_k\}$ is a set of edges in $G$ which contain a loop. Then 
$$\Psi^{I,J}_{G,K} =0 \quad  \hbox{ if } \quad  ( E\subset I\cup  K   \hbox{ or } E\subset J \cup K)  \quad \hbox{ and } \quad E \cap I \cap J= \emptyset \ .$$

\end{enumerate}

\subsection{Local structure} \label{sectlocal} We use these to deduce the local structure of $\Psi_G$ in some simple circumstances. Many more identities are derived in \cite{BrFeyn}. 
\begin{enumerate}
\item \emph{Local 2-valent vertex}. Suppose that $G$ contains a 2-valent vertex, whose neighbouring edges are labelled $1,2$. Then
$$\Psi^{12}_G=0\quad \hbox{ and }\quad \Psi_G^{1,2}=\Psi^1_{G,2} = \Psi^2_{G,1}$$
which imply that $\Psi_G= \Psi_{G\backslash 1\q 2} (\alpha_1+\alpha_2) + \Psi_{G\q \{1,2\}}$. 
In general, if $|I|=|J|$ are sets of edges such that $\{1,2\}\notin I\cup J\cup K$, then 
$$\Psi^{1I,2J}_{G,K} = \Psi^{I,J}_{G\backslash 1\q 2,K}=\Psi^{I,J}_{G\backslash 2\q 1,K}\ .$$
\item \emph{Doubled edge}. Suppose that $G$ contains doubled edges  $1,2$ (i.e., two edges which have the same set of  endpoints).  Then
$$\Psi_{G,12}=0 \quad \hbox{ and } \quad \Psi_G^{1,2}=\Psi^1_{G,2} = \Psi^2_{G,1}$$
which imply that $\Psi_G = \Psi_{G\backslash \{1,2\}} \alpha_1\alpha_2 + \Psi_{G\backslash 1\q 2} (\alpha_1+\alpha_2)$. 
In general, if $|I|=|J|$ are sets of edges such that $\{1,2\}\notin I\cup J\cup K$, then 
$$\Psi^{1I,2J}_{G,K} = \Psi^{I,J}_{G\backslash 1\q 2,K}=\Psi^{I,J}_{G\backslash 2\q 1,K}\ .$$
\item \emph{Local star}. Suppose that $G$ contains a three-valent vertex, whose neighbouring  edges  are labelled $1,2,3$. Then we have (\cite{BrFeyn}, Example 32)
$$\Psi^{123}_G=0\quad  \hbox{ and } \quad  \Psi^{12}_{G,3}= \Psi^{13}_{G,2} = \Psi^{23}_{G,1}$$
which follow from contraction-deletion. Furthermore, for  $\{a,b,c\}=\{1,2,3\}$  we  have the identities
$$\Psi^{ab,bc}_G= \Psi^{ab}_{G,c}=\ldots =\Psi^{bc}_{G,a} \quad \hbox{ and } \quad  \Psi^a_{G,bc} = \Psi^{a,c}_{G,b} +\Psi^{a,b}_{G,c}\ .$$
These identities propagate to higher order Dodgson polynomials. Let $i,j \notin \{1,2,3\}$.  Then for  all
$ \{a,b,c\}=\{a',b',c'\}=\{1,2,3\}$, we have (\cite{BrFeyn}, \S7.4):
$$\Psi^{abc, a ij}_G=0  \quad \hbox{and} \quad \Psi^{aci,bcj}_G =  \pm  \Psi^{i,j}_{G\backslash \{a',b'\}\q c'}\ .$$
\item \emph{Local triangle}.
Suppose that $G$ contains a triangle, with edges $1,2,3$. Then
$$\Psi_{G,123}=0\quad  \hbox{ and } \quad  \Psi^{1}_{G,23}= \Psi^{2}_{G,13} = \Psi^{3}_{G,12}$$
which follow from contraction-deletion. Furthermore, for  $\{a,b,c\}=\{1,2,3\}$  we  have the identities (\cite{BrFeyn}, Example 33):
$$\Psi^{a,b}_{G,c}= \Psi^{a}_{G,bc}=\ldots =\Psi^{b}_{G,ac} \quad \hbox{ and } \quad  \Psi^{ab}_{G,c} = \Psi^{ab,ac}_G +\Psi^{ab,bc}_G\ .$$
Now let $i,j \notin \{1,2,3\}$. For  all
$ \{a,b,c\}=\{a',b',c'\}=\{1,2,3\}$, we have
$$\Psi^{ab, ij}_{G,c}=0  \quad \hbox{and} \quad \Psi^{ai,bj}_{G,c} =  \pm  \Psi^{i,j}_{G\backslash  a' \q \{b',c'\}}\ .$$
\end{enumerate}

\subsection{The five-invariant}

\begin{defn} Let $i,j,k,l,m$ denote any five distinct edges in a graph $G$. The five-invariant of these edges, denoted ${}^5 \Psi_G(i,j,k,l,m)$ is defined to be the determinant
$${}^5\Psi_G(i,j,k,l,m) = \pm \det \left(
\begin{array}{cc}
 \Psi^{ij,kl}_{G,m} &   \Psi^{ik,jl}_{G,m}   \\
 \Psi^{ijm,klm}_{G} &   \Psi^{ikm,jlm}_{G}
\end{array}
\right)
$$
\end{defn}
It can be shown  that the five-invariant is well-defined, i.e., permuting the five indices $i,j,k,l,m$ only modifies the right-hand determinant by a sign.
In general, the 5-invariant is  irreducible of degree 2 in each Schwinger variable $\alpha_e$. However, in the case when three of the five edges $i,j,k,l,m$ form a star or a triangle,
it splits, i.e., factorizes into a product of Dodgson polynomials.

\begin{example} \label{example5split} Suppose that $G$ contains a triangle $a,b,c$. Then 
$${}^5\Psi_G(a,b,c,i,j) = \pm \det \left(
\begin{array}{cc}
 \Psi^{ab,ij}_{G,c} &   \Psi^{ai,bj}_{G,c}   \\
 \Psi^{abc,cij}_{G}  &\Psi^{aci,bcj}_{G}   
\end{array}
\right) =   \pm  \Psi^{i,j}_{G\backslash a\q \{b,c\}}  \Psi^{abc,cij}_{G} \ .
$$
It factorizes because $\Psi_{G,c}^{ab,ij}=0$ by the vanishing property for loops. By contraction-deletion,  
$\Psi^{ai,bj}_{G,c} = \Psi^{ai,bj}_{G\q c}$, and this is $\Psi^{i,j}_{G\backslash a\q \{b,c\}}$,   by the last equation of 
\S\ref{sectlocal}, (2), 
since  $a,b$   form a doubled edge in the quotient graph  $G\q c$.
\end{example}

\subsection{Denominator reduction} \label{sectdenomred}
Given a graph $G$ and an ordering $e_1, \ldots, e_{N_G}$ on its edges, we can extract a sequence of higher invariants, as follows.

\begin{defn}   Define $D^5_G(e_1,\ldots, e_5) = {}^5 \Psi_G(e_1,\ldots, e_5)$. 
 Let  $n\geq 5$ and suppose that we have defined $D^n_G(e_1,\ldots, e_n)$. Suppose furthermore that
$D^n_G(e_1,\ldots, e_n)$ factorizes into a product of linear factors in $\alpha_{n+1}$, i.e., it is of the form 
$(a\alpha_{n+1}+b) (c\alpha_{n+1}+d)$. Then we define
 $${D}^{n+1}_G(e_1,\ldots, e_{n+1}) =\pm( ad-bc)\ ,$$
 to be the resultant of the two factors of ${D}^n_G(e_1,\ldots, e_n)$. 
 A graph $G$ for which the polynomials ${D}^n_G(e_1,\ldots, e_{n})$ can be defined for all $n$ is called \emph{denominator-reducible}. 
 It can happen that $D^n_G(e_1,\ldots, e_n)$ vanishes. Then $G$ is said to have \emph{weight-drop}.
 \end{defn}

For general graphs above a certain loop order and   any ordering on their edges,  there will come a point where ${D}^n_G(e_1,\ldots, e_n)$ is irreducible (typically for $n=5$).
Thus the generic graph is not denominator reducible.  One can prove, as for the 5-invariant, that ${D}^n_G(e_1,\ldots, e_n)$ does not depend on the order of reduction of the variables,
although it may happen that the intermediate terms  ${D}^k_G(e_{i_1},\ldots, e_{i_k})$ may factorize for some choices of orderings and not others.

\section{The class of $X_G$ in the Grothendieck Ring of Varieties}

Let $k$ be a field. The Grothendieck ring of varieties $K_0(\mathrm{Var}_k)$ is  the free abelian group
generated by isomorphism classes  $[X]$, where $X$ is a separated  scheme of finite type over $k$,  modulo the inclusion-exclusion relation 
$[X]=[X\backslash Z] + [Z]$,  where $Z\subset X$ is a closed subscheme.  It has the  structure of a commutative ring   induced by the 
product relation $[X\times_k Y] = [X]\times [Y]$, with unit $1= [\mathrm{Spec}\, k]$. 
One defines the Lefschetz motive $\Lef$ to be the class of the affine line $[\A^1_k]$.

\begin{remark}
We only consider affine varieties here. If $f_1,\ldots, f_{\ell} \in k[\alpha_1,\ldots, \alpha_n]$ are polynomials, we denote by 
$[f_1,\ldots, f_{\ell}]$ the class in $K_0(\mathrm{Var}_k)$ of   the intersection of the hypersurfaces $V(f_1)\cap \ldots \cap V(f_{\ell})$ in affine space 
$\A^n_k$.  The dimension of the ambient affine space will usually be clear from the context.
\end{remark}

Let $G$ be a graph. Since the graph polynomial $\Psi_G$ (and, more generally, all Dodgson polynomials $\Psi^{I,J}_{G,K}$)  is defined over $\Z$, we can view the element
$[\Psi_G]$ in $K_0(\mathrm{Var}_k)$  for any field $k$. Most of the results below are  valid in this generality.
But at a certain point, we are obliged to  switch to point-counting functions since we require the use of the  Chevalley-Warning theorem (theorem \ref{thmCW}).
Recall that if $k$ is a finite field,   the point-counting map: 
\begin{eqnarray} \# : K_0(\mathrm{Var}_k) &  \rightarrow  &  \Z \nonumber \\
 {[}X] &\mapsto& \#X(k) \nonumber 
 \end{eqnarray}
 is well-defined, so results about the point-counts can be deduced from results in the Grothendieck ring, but not conversely (for example, it is not known if $\Lef$ is a zero-divisor).
 In this case, we shall denote by $[X]_q$ the point-counting function which associates to all prime powers $q$ the integers $ \# X(\F_q)$.

\subsection{Linear reductions}  The main observation of \cite{Stem} is that  the class in the Grothendieck ring of  polynomials which are linear in 
many of their variables    can be computed inductively by some simple  reductions.

\begin{lem}\label{lemlin}
 Let $f^1,f_1,g^1,g_1 \in k[\alpha_2,\ldots, \alpha_{n}]$ be  polynomials. Then

i).  $[f^1 \alpha_1+f_1] =  [f^1, f_1] \,\Lef+\Lef^{n-1} -[f^1]$

ii). $[f^1\alpha_1+f_1,g^1 \alpha_1+g_1] = [f^1, f_1, g^1, g_1] \, \Lef + [f^1g_1-f_1g^1]-[f^1, g^1]$
\end{lem}
Various proofs of this lemma  can be found in (\cite{SchnetzFq}, \cite{B-E-K} \S8, \cite{Stem} lemma 2.3, or $\S3.4$ of \cite{AM}).
Note that the quantity $ f^1g_1-g^1f_1$ is nothing other than the resultant $(\ref{resultantdef})$ with respect to $\alpha_1$ of  the polynomials $f^1\alpha_1+f_1 $ and $g^1\alpha_1 +g_1$.

Henceforth   let $G$ be connected.
  We call a graph simple if it has no vertices  of valency 
  $\leq 2$ (below left), multiple edges (below  
middle), or self-loops (below right).
\begin{center}
\fcolorbox{white}{white}{
  \begin{picture}(264,43) (138,-96)
    \SetWidth{1.0}
    \SetColor{Black}
    \Vertex(101,-79){2}
    \Vertex(191,-79){2}
    \Vertex(146,-79){2}
    \Line(101,-79)(191,-79)
    \Vertex(297,-79){2}
    \Vertex(252,-79){2}
    \Arc[clock](274.5,-106.625)(35.629,129.162,50.838)
    \Arc(274.5,-46.357)(39.646,-124.578,-55.422)
     \Vertex(367,-79){2}
      \Arc(377,-79)(10,-200,200)
    \Text(118,-71)[lb]{\Large{\Black{$e_1$}}}
    \Text(163,-71)[lb]{\Large{\Black{$e_2$}}}
    \Text(270,-67)[lb]{\Large{\Black{$e_1$}}}
    \Text(270,-98)[lb]{\Large{\Black{$e_2$}}}
  \Text(395,-82)[lb]{\Large{\Black{$e_1$}}}
  \end{picture}
}
\end{center}
\begin{lem} \label{SPlemma}
Let $G$ be a graph with a  subdivided edge $e_1,e_2$ (left). Then
\begin{equation} \label{Series} [\Psi_G]= \Lef [\Psi_{G\q e_1}]\ .
\end{equation}
Let $G$ be a graph with a doubled edge $e_1,e_2$
(middle). Then
\begin{equation}\label{Parallel}
[\Psi_G]= (\Lef-2) [\Psi_{G\backslash e_1}] +(\Lef-1) [\Psi_{G\backslash \{e_1,e_2\}}] + \Lef [ \Psi_{G\backslash e_1 \q e_2} ] + {\Lef}^{N_G-2}\ .
\end{equation}
Let $G$ be a graph with a self-loop $e_1$ (right). Then
\begin{equation} \label{Self_loop} [\Psi_G]= (\Lef-1) [\Psi_{G\backslash e_1}]+\Lef^{N_G-1}\ .
\end{equation}
\end{lem} 
\begin{proof}
These identities follow from the determination of the corresponding graph polynomials \S\ref{sectlocal} $(1), (2)$ and two applications of lemma \ref{lemlin} (see also \cite{AM}, \S4).
\end{proof}

An immediate consequence of the above lemma is that a smallest counter-example to Kontsevich's conjecture
will be a simple graph. Stronger restrictions are derived in \cite{Stem}.
Iterating the operations  $(\ref{Series})$ and $(\ref{Parallel})$ leads to explicit formulae for  $[X_G]$ as a polynomial in $\Lef$ when $G$ is a series-parallel graph  \cite{AM}.

\begin{propdef} \label{propdefc2}
  Let $G$ be a
graph such that $h_G \leq N_G-2$. Then  there exists an element  $c_2(G) \in K_0(\mathrm{Var}_k)/\Lef$  such that 
$$[\Psi_G] \equiv c_2(G) \Lef^2 \mod \Lef^3\ .$$

\end{propdef}

\begin{proof}   By induction we prove that for any  $F\in k[\alpha_1,\ldots, \alpha_n]$  of degree $<n$ which is linear in every variable $\alpha_i$, 
or any $G$ satisfying $h_G\leq N_G-2$ and any edge $e$ of $G$, there exist $a,b,c\in K_0(\mathrm{Var}_k)$ such that
\begin{enumerate}
\item $[F]\equiv a(F) \Lef \mod \Lef^2$
\item  $[\Psi_{G\backslash e}, \Psi_{G\q e} ]\equiv b(G,e) \Lef \mod \Lef^2$
\item $[\Psi_G] \equiv c(G) \Lef^2 \mod \Lef^3$.
\end{enumerate}
Proof of $(1)$.   The case where $n\leq 2$ is obvious. By linearity, let $F=f^1\alpha_1+f_1$. Lemma $\ref{lemlin}$ {\it(i)} implies that 
 $[F]\equiv [f^1,f_1] \Lef + \Lef^{n-1} - [f^1]$. 
Since  $f^1$ satisfies the required condition on the degrees, we can define $a$ inductively for $n>2$ by:
 $$a(F) = [f^1,f_1] - a(f^1)\ .$$
Proof of $(2)$.  From the contraction-deletion relations, we have $\Psi^1_G = \Psi^{12}_G \alpha_2 + \Psi^1_{G,2}$, and $\Psi_{G,1} = \Psi^{2}_{G,1} \alpha_2+ \Psi_{G,12}$.
The first Dodgson identity gives
$$\Psi^{1}_{G,2}\Psi^2_{G,1} - \Psi^{12}_G\Psi_{G,12} = (\Psi^{1,2}_G)^2.$$
Inserting this into  lemma $\ref{lemlin}$ {\it (ii)} gives
$$  [\Psi^1_G, \Psi_{G,1}] = [\Psi^{12}_G, \Psi^1_{G,2} , \Psi^2_{G,1}, \Psi_{G,12} ] \, \Lef +  [\Psi^{1,2}_G] - [\Psi^{12}_{G}, \Psi^2_{G,1}] \ . $$
We have $\deg \Psi^{1,2}_G = h_G -1 < N_G-2$ so $(1)$ applies to $\Psi^{1,2}_G$. If $G$ is not connected, define $b(G,1)$ to be $0$. Otherwise, define $b$ inductively by:
$$b(G,1)= a(\Psi^{1,2}_G) - b(G\backslash 2, 1) + [\Psi^{12}_G, \Psi^1_{G,2} , \Psi^2_{G,1}, \Psi_{G,12} ]\ .$$
If $G\backslash 2$ is  connected,  Euler's formula shows that  $G\backslash 2$ satisfies the required condition on the degree. The initial case when $h_G=0$, i.e.\
$G$ is a tree, is obvious.

Proof of $(3)$. By contraction-deletion we write $\Psi_G = \Psi^1_G \alpha_1 + \Psi_{G,1}$. By lemma $\ref{lemlin}$ {\it (i)}, 
$[\Psi_G] = [\Psi^1_G, \Psi_{G,1}] \Lef - [\Psi^1_G] + \Lef^{N_G-1}$, so for $N_G> 2$ we  inductively define 
$$c(G) = b(G,1) - c(G\backslash 1)$$
if $G$ is connected, and set $c(G)=0$ otherwise.  The case $N_G= 2$ is obvious.
\end{proof}
Note that in the opposite case,  if $G$ is connected and satisfies $h_G> N_G-2$, then $G$ has at most two vertices and is essentially uninteresting.
 
\begin{cor} \label{cor2valent} Suppose that $G$ has a 2-valent vertex
and $h_G\leq N_G-3$. Then $c_2(G) \equiv 0 \mod \Lef$.
\end{cor}
\begin{proof} Using  lemma \ref{SPlemma}, we can write $[\Psi_G] \equiv \Lef [\Psi_{G  \q e}] \equiv 0 \mod \Lef^3$
since we have  $h_{G\q e}=h_G\leq N_{G\q e}-2$. 
\end{proof}

\begin{remark}
Below  we give a formula for $c_2(G)$  when $2h_G\leq N_G$. The quantity $2h_G-N_G$ is connected to the physical `superficial degree
of divergence' in space-time dimension 4. Graphs with $2h_G<N_G$ are superficially convergent. The physically interesting case
is superficial log-divergence $2h_G=N_G$. Primitive-divergent graphs belong to this class.
\end{remark} 
\subsection{Three-valent vertices} \label{sect3vv}
Our approach to studying $[\Psi_G]$ uses the existence of a
vertex with low degree. Note that whenever
\begin{equation}\label{3vexists} 2h_G-2<N_G
\end{equation} the minimum vertex-degree is $\leq 3$.
To see this, note that   Euler's formula for a connected graph  implies that $N_G-V_G=h_G-1$. If $\alpha$ denotes the average degree of the vertices of $G$,
then $N_G=\frac{\alpha}{2} V_G$, and 
 $(\ref{3vexists})$ implies that $\alpha<4$.
\begin{figure}[ht!]
  \begin{center}
    \leavevmode
    \epsfxsize=3.0cm \epsfbox{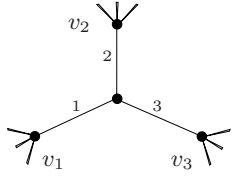}
 \put(-30,21){{\tiny $3$}}
 \put(-64,21){ {\tiny $1$}}
 \put(-75,0){ {\small $v_1$}}
 \put(-49,40){{\tiny  $2$}}
 \put(-62,52){{\small  $v_2$}}
  \put(-23,0){{\small $v_3$}}
  \end{center}
  \caption{A three-valent vertex}
\end{figure}

The case of a two-valent vertex was dealt with in \S\ref{sectlocal}. 
The case of a three-valent vertex is more complicated but  still implies that $\Psi_G$ has a simple structure.

\begin{defn}\label{defnfi}
Let $v_1,v_2,v_3$ be any three vertices in $G$ which form a  three-valent vertex as shown above. Following \cite{BrFeyn},  we will use the notation:
$$f_0 = \Psi_{G \backslash \{1,2\} \q 3} \ , \ f_1 =\Psi^{2,3}_{G,1} \ , \ f_2 = \Psi^{1,3}_{G,2}  \ , f_3 = \Psi^{1,2}_{G,3} \  , \ f_{123}=\Psi_{G\q \{1,2,3\} }  \  .$$ 
\end{defn}

\begin{lem} \label{lemvtxstructure}
In this case,  the graph polynomial of $G$ has the following structure:
$$\Psi_G=f_0 (\alpha_1\alpha_2+\alpha_1\alpha_3+\alpha_2\alpha_3)+(f_1+f_2)\alpha_3+(f_1+f_3)\alpha_2+(f_2+f_3)\alpha_1+ f_{123}$$
where the polynomials $f_i$ satisfy the equation
\begin{equation} \label{fids}
f_0f_{123}= f_1f_2+f_1f_3+f_2f_3\ .
\end{equation}
\end{lem}
\begin{proof} The general shape of the polynomial comes from the contraction-deletion relations, and  \S\ref{sectlocal} $(3)$ (or ex. 32 in \cite{BrFeyn}).
Equation $(\ref{fids})$ is merely a restatement of the first Dodgson identity for $G\q 3$ which gives
$(\Psi^{1,2}_{G,3})^2 = \Psi^1_{G,23}\Psi^2_{G,13}-\Psi^{12}_{G,3}\Psi_{G,123}.$
Using the definitions of $f_i$ this translates as 
\begin{equation}\label{2}
f_3^2 = (f_2+f_3)(f_1+f_3)-f_0f_{123}\ .
\end{equation}
\end{proof}

\begin{prop} \label{prop3valentvtx} Suppose that $G$ contains a three-valent vertex, and let $f_i$ be given by definition \ref{defnfi}.
Then 
$$ [\Psi_{G}]  = \Lef^{N_G-1} + \Lef^3 [f_0,f_1,f_2,f_3,f_{123}]  -  \Lef^2[f_0,f_1,f_2,f_3]  $$
\end{prop}

\begin{proof} Let $\beta_i = f_0 \alpha_i + f_i$, for $i=1,2,3$. It follows from $(\ref{fids})$ that
$$ f_0\Psi_G =\beta_1\beta_2+\beta_1\beta_3+\beta_2\beta_3 \ .$$
The right-hand side is the graph polynomial of the sunset graph (a triple edge). It defines a 
quadric in $\A^3$ whose class is $\Lef^2$. 
It follows that if $U$, and $U'$ denote the open set $f_0\neq 0$ in $\A^{N_G}$, and in $\A^{N_G-3}$ resp., we have $[X_G\cap U]=\Lef^{2}[U']$.
On the complement  $V(f_0)$, the graph polynomial $\Psi_G$ reduces to the equation
$$(f_1+f_2)\alpha_3+(f_1+f_3)\alpha_2+(f_2+f_3)\alpha_1+ f_{123}$$
which defines  a family of  hyperplanes in $\A^3$.  Thus, consider the fiber of the projection $X_G\cap V(f_0) \rightarrow \A^{N_G-3}\cap V(f_0)$.
In the generic case this is a hyperplane whose class is $\Lef^2$. Otherwise, $f_1,f_2,f_3$ vanish and
there are two possibilities: if $f_{123}=0$ the fiber is isomorphic to  $\A^3$, otherwise it is empty.  We have
$$[X_G\cap V(f_0)] = \Lef^3 [f_0,f_1,f_2,f_3,f_{123}] +\Lef^2([f_0] - [f_0,f_1,f_2,f_3])$$
Writing $[X_G]=[X_G\cap U ] + [X_G\cap V(f_0)]$ and $[U']=\Lef^{N_G-3}-[f_0]$ gives the result.
\end{proof} 

In particular, if $G$ has a three-valent vertex and $N_G\geq 4$ then 
\begin{equation} \label{c2symmetricform} 
c_2(G) \equiv -[f_0,f_1,f_2,f_3] \mod \Lef\ . 
\end{equation}

\begin{lem}\label{lem2}
Let $G$  satisfy $h_G+3\leq N_G$, where $N_G\geq 4$, and contain  a three-valent vertex whose neighbouring edges are numbered $1,2,3$.  Then 
\begin{equation} \label{c2at3stage} c_2(G) \equiv [\Psi^{1,2}_{G,3}, \Psi^{13,23}_{G}]  \mod \Lef\ . \end{equation}
\end{lem}
\begin{proof} 
We use the explicit expression for $\Psi_G$  in lemma $\ref{lemvtxstructure}$ and  the relations in \S\ref{sectlocal} $(3)$.
It follows from $(\ref{fids})$ and inclusion-exclusion that:
  $$[f_0,f_3] =[f_0,f_1f_2,f_3] =  [f_0,f_1,f_3]+[f_0,f_2,f_3] - [f_0,f_1,f_2,f_3]$$
On the other hand,  $[f_0,f_1+f_3] = [f_0,f_1+f_3,f_3^2]$ by equation $(\ref{2})$, and  so we have  $[f_0,f_1+f_3]=[f_0,f_1,f_3]$.
By  contraction-deletion, we can write
 $$[f_0,f_1+f_3] = [\Psi^{12}_{G,3}, \Psi^{2}_{G,13}] =[\Psi^1_{G'}, \Psi_{G',1}]\ ,$$
where $G'=G\backslash 2 \q 3$. Either $G'$ is not connected, or else $h_{G'} \leq N_{G'}-2$ by the  assumption on the loop number of $G$,
and so  the previous expression  vanishes modulo $\Lef$ by statement $(2)$ in the proof of proposition-definition $\ref{propdefc2}$.
The same is true for $[f_0,f_1+f_2]$ by symmetry.  We have therefore shown that 
  $$-[f_0,f_1,f_2,f_3] \equiv  [f_0,f_3]  \mod \Lef\ .$$
The lemma follows from  $(\ref{c2symmetricform})$ and   $\Psi^{1,2}_{G,3} =f_3 $ and $\Psi^{13,23}_{G} =\Psi^{12}_{G,3} =f_0 $.
\end{proof}

\subsection{Counting points over finite fields} 
For any prime power $q$, let $\F_q$ denote the  field with $q$ elements. Given polynomials $P_1,\ldots, P_{\ell}\in \Z[\alpha_1,\ldots, \alpha_n]$, let  
$$[P_1,\ldots, P_{\ell}]_q \in \mathbb{N}\cup \{0\}$$
denote the number of points on the affine variety $V(\overline{P_1},\ldots, \overline{P_{\ell}}) \subset \F_q^{n}$, where $\overline{P}_i$ denotes the reduction of $P_i$ modulo $p$ (the characteristic of $\F_q)$.
Recall  the Chevalley-Warning theorem (e.g., \cite{Se})  on the point-counts of polynomials of small degrees. 

\begin{thm} \label{thmCW} Let $P_1,\ldots, P_{\ell} \in \Z[\alpha_1,\ldots, \alpha_n]$ such that $\sum_{i=1}^{\ell} \deg P_i <n$. Then 
$$[P_1,\ldots, P_{\ell}]_q \equiv  0 \mod q\ .$$
\end{thm}
It is natural to ask if  the Chevalley-Warning theorem lifts to the Grothendieck ring of varieties. We were unable to find such a result in the literature.

\begin{ques} \label{conjCW} For which fields $k$ is the following statement true: Let $P_1,\ldots, P_{\ell}$ be polynomials satisfying the above condition on their degrees.
Then $[V(P_1,\ldots, P_{\ell})] \equiv 0 \mod \Lef$ in $K_0(\hbox{Var}_k)$?
\end{ques}
In an earlier version of this paper we cautiously conjectured this to be true for all $C_1$ fields  (see the examples in   \cite{EV}), 
which  (as pointed out to us by a referee) would imply the result  for all fields of finite characteristic.
Lacking strong evidence for this, it is perhaps more prudent to assume $k$ to be algebraically closed.
In any case,  since a  geometric Chevalley-Warning theorem is unavailable,
we henceforth work with  point-counting functions rather than with elements in the Grothendieck ring of varieties.
It turns out that for many of the results below, one can in fact circumvent this question by elementary arguments.
Nevertheless, we now set
$$c_2(G)_q = [\Psi^{13,23}_G,  \Psi^{1,2}_{G,3}]_q \mod q$$
viewed as a map from all prime powers $q$ to $\Z/q \Z$, and where $1,2,3$ forms a three-valent vertex as above. Below we show that the formula remains valid  for any
set of three edges $1,2,3$. We have $[\Psi_G]_q \equiv c_2(G)_q q^2 \mod q^3$.

\begin{lem}\label{lem1} Suppose that $f=f^1
\alpha_1+f_1$ and $g=g^1\alpha_1+g_1$ are polynomials in $\Z[\alpha_1,\ldots, \alpha_n]$ such that $\deg f+ \deg g 
\leq n$, which are linear in a variable $\alpha_1$. Then
$$[f,g]_q \equiv [f^1g_1-f_1g^1]_q \mod q\ .$$
If the resultant  has a non-trivial factorization $f^1g_1-f_1g^1=ab$,  then 
$$[f,g]_q \equiv- [a,b]_q \mod q\ .$$
\end{lem}
\begin{proof} By lemma \ref{lemlin} {\it (ii)},  $[f,g]_q= q [f^1,f_1,g^1,g_1]_q +[f^1g_1-f_1g^1]_q-[f^1,g^1]_q$. Since $f^1,g^1 \in \Z[\alpha_2,\ldots, \alpha_n]$ have total
degree $\deg(f^1)+\deg(g^1)\leq n-2$, this is congruent to $[f^1g_1-f_1g^1]_q$ mod $q$ by  theorem $\ref{thmCW}$.

By inclusion-exclusion we have
$[ab]_q=[a]_q+[b]_q-[a,b]_q$. The  factorization  is non-trivial  if and only if $\deg a $ and $\deg b$ are strictly smaller than $\deg ab$. By theorem $\ref{thmCW}$,
this implies that   $[a]_q$ and $[b]_q$ vanish mod $q$, giving the second statement. 
\end{proof}

\begin{cor}\label{cor1}
Let $G$ be a connected graph such that $2h_G\leq N_G$, $N_G\geq5$,  and  let $1,2,3$ be any distinct edges of $G$. Then
$$[\Psi^{1,2}_{G,3},\Psi^{13,23}_G]_q\equiv-[{}^5\Psi_G]_q\mod q\ ,$$
where the 5-invariant is taken with respect to any set of five edges of $G$. In particular, the point-counts of all 5-invariants
are equivalent mod $q$.
\end{cor}
\begin{proof}
First assume that the edges 1,2,3 are a subset of the edges in the 5-invariant ${}^5\Psi(1,2,3,4,5)$.
We have  $\deg(\Psi^{1,2}_{G,3}) = h_G-1$  and  $ \deg \Psi^{13,23}_{G} =h_G-2$
giving total degree $2h_G-3$, whereas the ambient affine space has dimension $N_G-3$.
Applying the previous lemma to equation $(\ref{c2at3stage})$  gives 
$$ [\Psi^{1,2}_{G,3}, \Psi^{13,23}_G]_q \equiv -[\Psi_G^{13,24},\Psi_G^{14,23}]_q \mod q$$
by the Dodgson identities.  Applying the previous lemma one more time gives
$$-[\Psi_G^{13,24},\Psi_G^{14,23}]_q  \equiv -[{}^5\Psi_{G}(1,2,3,4,5)]_q \mod q$$
by definition of the five-invariant as a resultant.
Since the edges 4 and 5 are arbitrary we see that the point-counts of 5-invariants are equivalent mod $q$ whenever
they have three edges in common. By considering chains of overlapping edge-sets the same is true for any 5-invariants.
\end{proof}
In particular, if $G$ has a three-valent vertex, then $c_2(G)_q \equiv [\Psi^{13,23}_G,\Psi^{1,2}_{G,3}]_q \mod q$ for \emph{any} three edges 1,2,3 of $G$ which do not necessarily
meet  the three-valent vertex.
\begin{thm} \label{corDenom}
Let $G$ be a connected graph with  $2h_G\leq N_G$, $N_G\geq5$. Suppose that $D^n_G(e_1,\ldots, e_n)$ is the result of
the denominator reduction after $n<N_G$ steps. Then
\begin{equation}\label{3}
c_2(G)_q \equiv  (-1)^n [D^n_G(e_1,\ldots, e_n)]_q \mod q\ .
\end{equation}
If $G$ has weight drop or $2h_{G} < N_G\geq4$, then $c_2(G)_q \equiv 0 \mod q$.
\end{thm}
\begin{proof}
Suppose first that $G$ has a three-valent vertex.  Equation (\ref{3})  follows by induction from corollary \ref{cor1} by
 applying the denominator reduction \S\ref{sectdenomred}.
 There are two cases to consider: if the factorization in the denominator reduction is non-trivial
and $n<N_G$ then the induction step follows from  lemma \ref{lem1}. If the factorization in the denominator reduction is trivial (the denominator is of degree one
in the reduction variable), then it follows from lemma \ref{lemlin} {\it i}).

If $G$ does not have a three-valent vertex, then 
since $2h_G\leq N_G$,  one of the following situations must hold (see the argument
in \S\ref{sect3vv}): (i) $G$ has a two-valent vertex, (ii) $G$ has a one-valent vertex, or (iii) $G$ has a self-loop
connected to a single edge (forming a degenerate three-valent vertex).

In the case (i), we conclude from $2h_G\leq N_G$ and $N_G\geq 5$  that $h_G\leq N_G-3$  and so $c_2(G)_q\equiv0\mod q$ by corollary
\ref{cor2valent}. Likewise,  the 5-invariant vanishes if it contains a two-valent vertex
(see \cite{BrFeyn} lemma 92). The same argument holds trivially in the cases (ii) and (iii).

In case of a weight drop the right hand side of (\ref{3}) vanishes, hence $c_2(G)_q\equiv0\mod q$. The 5-invariant is
of degree $2h_G-5$ in $\A^{N_G-5}$. If $2h_G<N_G\geq5$ we have $c_2(G)_q\equiv0\mod q$ by theorem \ref{thmCW}.
If $2h_G<N_G=4$ then $h_G\leq 1$. The  hypersurface $X_G$ is either empty or a hyperplane in $\A^4$,
hence $c_2(G)_q\equiv0\mod q$. 
\end{proof}
Notice that the proofs rely on the fact that  the terms  $D^n_G(e_1,\ldots, e_n)$ in the denominator reduction are of degree exactly equal to the dimension  of the ambient space,
and therefore lie on the limit of the Chevalley-Warning theorem (the Calabi-Yau condition for the  associated projective varieties).

\section{Properties of the $c_2$-invariant}
We state some known  and conjectural properties of the $c_2$-invariant of a graph.
Throughout  this section $G$ is a graph with $2h_G\leq N_G$ and at least five edges.

\subsection{Triviality of $c_2(G)$} The following results follow from theorem \ref{corDenom}.
\begin{lem} If $G$  as above   has a
doubled edge then $c_2(G)_q \equiv 0 \mod q$. 
\end{lem}
\begin{proof} If $G$ has a doubled edge $e_1,e_2$, then any five-invariant ${}^5\Psi_G(i_1,\ldots, i_5)$ 
where $e_1,e_2\in \{i_1,\ldots, i_5\}$ necessarily vanishes (\cite{BrFeyn} lemma 90).
\end{proof} 
Recall that $G$ is called 2-vertex reducible if there is a pair of distinct  vertices  such that removing them (and their incident edges) causes the graph  to disconnect.
\begin{prop} Let $G$ as above be  2-vertex reducible. Then $c_2(G)_q \equiv 0 \mod q$.
\end{prop}
\begin{proof} It is proved in \cite{WD}, proposition 36, that such a graph has weight drop.
\end{proof}
\begin{prop} \label{propdenomcase} 
If $G$ is  denominator reducible, and non-weight drop, then $c_2(G)_q$ $\equiv (-1)^{N_G-1}\mod q$ except for finitely many primes $p$ for which $c_2(G)_{p^n}\equiv0 \mod p^n$.
\end{prop} 
\begin{proof}
If $G$ is non-weight drop denominator reducible then there exists a degree one homogeneous polynomial
$D_G^{N_G-1}(e_1,\ldots,e_{N_G-1})=c\alpha_{N_G}$ with $0\neq c\in \Z$. For primes $p|c$ we have $[c\alpha_{N_G}]_{p^n}=p^n$,
otherwise $[c\alpha_{N_G}]_q=1$. The result follows from theorem \ref{corDenom}.
\end{proof}
\subsection{Double triangle reduction}
Consider a graph $G$ which contains seven edges $e_1,\ldots, e_7$ arranged in the configuration shown below on the left (where anything may be attached to vertices $A$-$D$).
The double triangle reduction of $G$ is the graph $G'$ obtained by replacing these seven edges with the configuration of five edges $e_1',\ldots, e_5'$ as shown below on the right.
The following theorem was proved in \cite{WD}.
\vspace{-0.05in}
\begin{center}
\fcolorbox{white}{white}{
  \begin{picture}(302,120) (6,-5)
    \SetWidth{1.0}
    \SetColor{Black}
    \Line(44,52)(66,95)
    \Line(65,95)(87,52)
    \Line(44,52)(66,9)
    \Line(65,9)(87,52)
    \Line(16,52)(116,52)
    \Line(222,95)(251,52)
    \Line(251,52)(222,9)
    \Line(230,52)(273,52)
    \Line(222,95)(222,9)
    \Line(187,52)(216,52)
    \Vertex(44,52){3}
    \Vertex(87,52){3}
    \Vertex(251,52){3}
    \Arc(16,52)(3,270,630)
    \Arc(65.5,95)(3,270,630)
    \Arc(65.5,9)(3,270,630)
    \Arc(116,52)(3,270,630)
    \Arc(187,52)(3,270,630)
    \Arc(222,95)(3,270,630)
    \Arc(273,52)(3,270,630)
    \Arc(222,9)(3,270,630)
    \Text(20,10)[lb]{\Large{\Black{$G$}}}
    \Text(180,10)[lb]{\Large{\Black{$G'$}}}
    \Text(16,58)[lb]{{\Black{$A$}}}
    \Text(65,102)[lb]{{\Black{$B$}}}
    \Text(65,-6)[lb]{{\Black{$C$}}}
    \Text(116,38)[lb]{{\Black{$D$}}}
    \Text(180,58)[lb]{{\Black{$A$}}}
    \Text(273,38)[lb]{{\Black{$D$}}}
    \Text(222,102)[lb]{{\Black{$B$}}}
    \Text(222,-6)[lb]{{\Black{$C$}}}
  \end{picture}
}
\end{center}

\begin{thm}  Let $G'$  be a double triangle reduction of $G$. Then 
$$D^7_G (e_1,\ldots, e_7) = \pm D^5_{G'} (e'_1,\ldots, e'_5)\ .$$
\end{thm}
\begin{cor} Let $G, G'$  be as above,  with $2h_G\leq N_G$. Then
$$c_2(G)_q \equiv c_2(G')_q \mod q\ .$$
\end{cor} 
Since the double-triangle reduction violates planarity, this is the first hint that the genus of a graph is \emph{not} the right invariant for understanding its periods.
\subsection{The completion relation}
It follows from a simple application of Euler's formula that a primitive-divergent graph $G$ in $\phi^4$ with more than
one loop has exactly four three-valent vertices $v_1,\ldots, v_4$, and all remaining vertices
have valency 4.  The completion of $G$ is defined to be the graph $\widehat{G}$ obtained by adding a new vertex  $v$ to $G$ and connecting it to $v_1,\ldots, v_4$ \cite{SchnetzCensus}.
The resulting graph is 4-regular. 

\begin{conj} \label{completionconjecture} Let $G_1,G_2$ be two primitive-divergent graphs in $\phi^4$ and suppose that $\widehat{G_1} \cong \widehat{G_2}$.  Then
$c_2(G_1) \equiv  c_2(G_2) \mod \Lef$\ .
\end{conj} 
The motivation for this conjecture comes from the result \cite{SchnetzCensus} that the corresponding residues are the same: $I_{G_1}=I_{G_2}$. 
 Once again, the completion relation does not respect the genus of a graph.

\section{Mixed Tate families: Graphs of Vertex-width 3}

When $G$ contains  sufficiently many triangles and three-valent vertices, we show that $[\Psi_G] \in K_0(\mathrm{Var}_k)$ is a polynomial in $\Lef$ which can be computed 
inductively.

\subsection{The vertex-width of a graph}\label{sectvw} Throughout, $G$ is a connected graph.

\begin{defn}   Let $\mathcal{O}$ be an ordering on the edges of $G$. It gives rise to a filtration
$$\emptyset=G_0 \subset G_1 \subset \ldots \subset G_{N-1} \subset G_N = G$$
of
 subgraphs of $G$, where $G_i$ has exactly $i$ edges. To any such filtration we obtain  a sequence of integers
$v^{\mathcal{O}}_i = \hbox{number of vertices of } G_i \cap (G\backslash G_i).$
We say that $G$ has \emph{vertex-width}  at most  $n$ if there exists an ordering $\mathcal{O}$ such that $v^{\mathcal{O}}_i \leq n$ for all $i$ \cite{BrFeyn}.
\end{defn}

For example, a row of boxes with vertices $a_1,\ldots, a_n, b_1,\ldots, b_n$ and edges $\{a_i,b_i\}$,  $\{a_i,a_{i+1}\}$,  $\{b_i,b_{i+1}\}$,   has vertex-width two.
The wheels and zig-zag graphs (below) have vertex-width $\leq 3$. Bounding the vertex-width is a strong constraint on a graph, and one can show that the set of planar graphs
have arbitrarily high vertex-width. In \cite{BrFeyn}  it was shown  that the relative cohomology of the graph hypersurface for   graphs of vertex width $\leq 3$  is mixed Tate,
and  that the periods are multiple polylogarithms. Here we explain how to compute the class of $[X_G]$ as a polynomial in $\Lef$ for such graphs.
For this, it is not enough to consider  only the classes $[\Psi_H]\in K_0(\mathrm{Var}_k)$, where $H$ are minors of $G$, and we are forced to introduce a new invariant:
\begin{defn}
Let $e_1,e_2,e_3$ be any three edges in $G$ which form a  three-valent vertex. If $f_0,f_1,f_2,f_3,f_{123}$ are given  by definition \ref{defnfi}, we set
\begin{equation} \langle G \rangle_{e_1,e_2,e_3}= [f_0,f_1,f_2,f_3,f_{123}]\end{equation}
in $\A^{N_G-3}$. Sometimes we shall write $\langle G\rangle_v$ if $v$ is the 3-valent vertex meeting edges $e_1,e_2,e_3$.
\end{defn}
We first consider recurrence relations for $[\Psi_G]$ (which also involve invariants $\langle H\rangle$ for minors $H$ of $G$), and then recurrence relations for $\langle G\rangle$
(which also involve invariants $[\Psi_H]$ for minors $H$ of $G$).

\subsection{Reduction of $[\Psi_G]$}
The  two main cases are split triangles and split vertices. 

\subsubsection{Split Vertices}   Let $G$ be any graph containing a three-valent vertex  (left),  and let $G'$ be the graph obtained by splitting that vertex in two (right).
An empty (white) vertex indicates that there can be other  edges connected to it which are not  drawn on the diagram (anything can be attached to $v_1,v_2,v_3$).
\begin{figure}[ht!]
 \begin{center}
    \leavevmode
    \epsfxsize=10.0cm \epsfbox{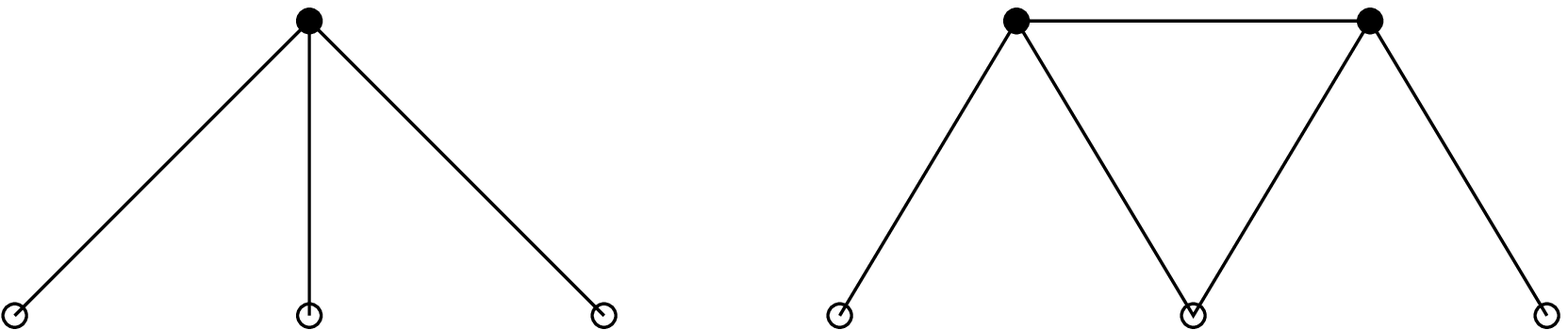}
\put(-266,30){$1$} \put(-237,30){$2$} \put(-199,30){$3$}
 \put(-290,-10){$v_1$} \put(-235,-10){$v_2$} \put(-180,-10){$v_3$}
 \put(-140,-10){$v_1$} \put(-75,-10){$v_2$} \put(-15,-10){$v_3$}
  \put(-125,30){$4$} \put(-82,30){$1$} \put(-61,30){$3$}\put(-17,30){$5$} \put(-71,47){$2$} \put(-100,62){$v_4$}\put(-50,62){$v_5$}
  \put(-300,50){\large{$G$}}  \put(-150,50){\large{$G'$}}
  \end{center}
  \caption{}
 \end{figure}

\begin{thm}  \label{thmvertexsplit} The class of the graph polynomial of $G'$ can be written explicitly in terms of the invariant $\langle G \rangle_{1,2,3} $, and the classes
of minors of $G$:
$$[\Psi_{G'} ] + (\Lef-\Lef^2)\big([\Psi_{G,2}] - [\Psi^{13}_{G,2}]\big)+ (\Lef-1) [\Psi_{G}] = (\Lef^5-\Lef^4)\langle G \rangle_{1,2,3}  +\Lef^{N_G-2}(\Lef^3+\Lef-1)$$
\end{thm}
\begin{proof} The structure of the graph polynomial of $G'$ can be obtained as follows. Since $v_4$ is a three-valent vertex in $G'$, it follows that $\Psi_{G'}$ must be of the shape
given in lemma $\ref{lemvtxstructure}$, for some polynomials $f'_0, f'_1,f'_2,f'_4,f'_{124}$ relative to the edges $1,2,4$. By contraction-deletion relations, one easily sees that
$$f'_0= f_0 (\alpha_3+\alpha_5) + (f_2+f_3) \quad , \quad f'_{124} = f_{123} \alpha_3 + (f_1 +f_2) \alpha_3\alpha_5 $$
$$f'_1=f_2\alpha_3\quad , \quad f'_2=f_3\alpha_3 +f_0 \alpha_3\alpha_5\quad , \quad f'_4=f_{123} + f_1 \alpha_3 +(f_1+f_2)\alpha_5 \ ,$$
where $f_0,f_1,f_2,f_3,f_{123}$ satisfy $(\ref{fids})$ and are the invariants of the three-valent vertex formed by edges $1,2,3$ of $G$.
By proposition $\ref{prop3valentvtx}$ we know that  $[\Psi_{G'}]$ is given by  
$\Lef^{N_{G'}-1}+ \Lef^3 [f'_0,f'_1,f'_2,f'_4,f'_{124}] - \Lef^2 [f'_0,f'_1,f'_2,f'_4]$. The conclusion of the theorem follows by a brute force calculation by exploiting the
inclusion-exclusion relations, identity $(\ref{fids})$, and  reducing out the linear variables $\alpha_3,\alpha_5$ using lemma \ref{lemlin} {\it (ii)}. 
\end{proof}
Any inductive procedure to compute the class of a split-vertex graph $G'$  is blocked by the presence of an invariant $\langle G\rangle$.
However, modulo $\Lef^4$  it drops out.
\begin{cor} Suppose that $
N_G\geq 6$. Then $c_2(G')  \equiv  c_2(G) \mod \Lef$. If  $[\Psi_G]$ is of the form
$[\Psi_G] \equiv c_3(G) \Lef^3 + c_2(G) \Lef^2 \mod \Lef^4$, then so is $[\Psi_G']$ and we have:
$$c_3(G')-c_3(G)  \equiv    c_2(G\backslash \{1,3\}\q 2) -     c_2(G \q 2) -c_2(G)  \mod \Lef  $$
\end{cor}
\begin{proof} This follows from theorem \ref{thmvertexsplit} and proposition-definition \ref{propdefc2}.
\end{proof}
Iterating this corollary leads, for example,  to an inductive way to compute the coefficient  $c_3$ of $\Lef^3$ for certain classes of graphs which are polynomials in $\Lef$.

\subsubsection{Split triangles} Let $G'$ be a graph of the shape depicted below (right), and let $G$ denote the subgraph obtained by deleting edges 4 and 5.

\begin{figure}[ht!]
 \begin{center}
    \leavevmode
    \epsfxsize=10.0cm \epsfbox{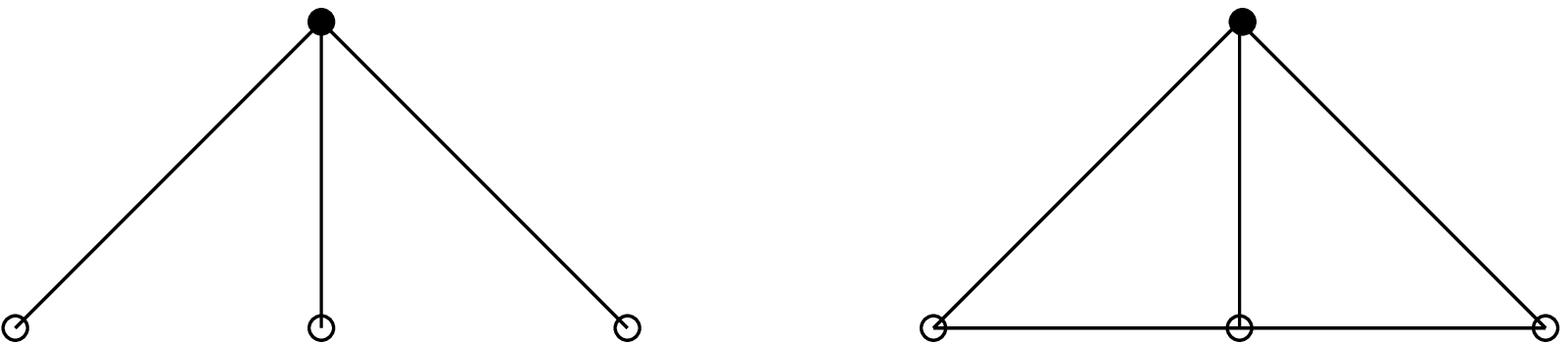}
\put(-266,30){$1$} \put(-233,30){$2$} \put(-196,30){$3$}
 \put(-287,-10){$v_1$} \put(-232,-10){$v_2$} \put(-175,-10){$v_3$}
 \put(-120,-10){$v_1$} \put(-65,-10){$v_2$} \put(-5,-10){$v_3$}
  \put(-90,7){$4$} \put(-95,30){$1$} \put(-67,30){$2$}\put(-27,30){$3$} \put(-37,7){$5$}
  \put(-285,45){\large{$G$}}  \put(-125,45){\large{$G'$}}
  \end{center}
  \caption{}
 \end{figure}

\begin{thm}\label{thmsplittriangle} 
 Let $G'$, $G$ be as above, and let $H=G'\backslash \{1,3\}\q 2$, and $\triangle= G'\backslash 2 \q  3$. The following equation relates $[\Psi_{G'}]$ to 
the classes of minors of $G'$, and  $\langle G\rangle_{1,2,3}$: 
$$[\Psi_{G'}] +[\Psi^4_{G'}]+[\Psi^5_{G' }] +  [\Psi^{45}_{G'}] + \Lef \big( [\Psi_{H}] + [\Psi^{4}_{H}] + [\Psi^{5}_{H}]  + [\Psi^{45}_{H}]+ 
[\Psi^{4}_{\triangle,1}] +[\Psi^{45}_{\triangle,1}]
\big) $$
$$- (\Lef^3-\Lef^2)\big( [\Psi_{H,4}]  +[\Psi_{H,5}]+ [\Psi^{4}_{H,5}] +[\Psi^{5}_{H,4}]  + [\Psi_{H,45}]  \big) - \sum_{T\subseteq \{1,4,5\}} [\Psi^T_{{\triangle}}]  $$
$$= (\Lef^5-\Lef^4) \langle G \rangle_{1,2,3}  +(\Lef^4+3\, \Lef^2 -\Lef-1) \, \Lef^{N_H-2}$$
where the sum  is over all $8$ subgraphs $T$ of ${\triangle}$ obtained by deleting the  edges $1,4,5$. 
\end{thm}

\begin{proof}  We omit the proof, which is similar to the proof of theorem \ref{thmvertexsplit}.
\end{proof} 

\subsection{Recurrence relations for $\langle G\rangle$} \label{sectGnewinv}
 It turns out that  $\langle G\rangle$ satisfies  recurrence relations with respect to a set of four edges. Consider the three graphs below:

\begin{figure}[ht!]
 \begin{center}
    \leavevmode
    \epsfxsize=12.0cm \epsfbox{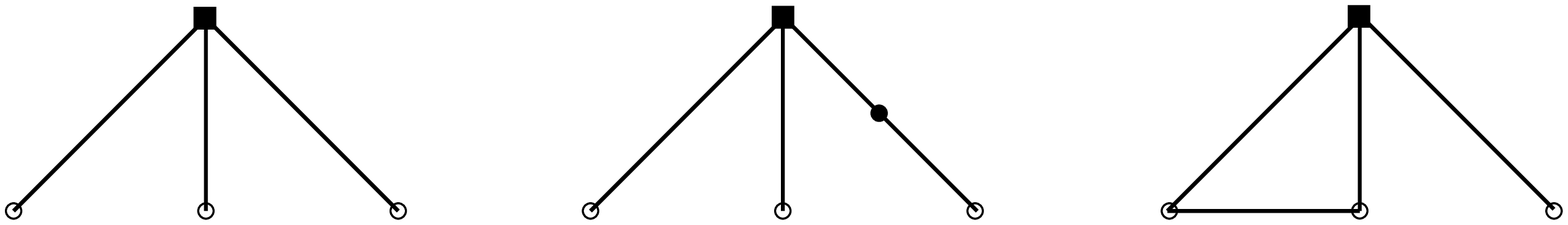}
   \put(-303,-15){$G$}  
\put(-325,25){$1$} \put(-305,25){$2$} \put(-275,25){$3$}
   \put(-177,-15){$G_s$}  
\put(-200,25){$1$} \put(-180,25){$2$} \put(-157,33){$3$} \put(-137,15){$4$}
     \put(-50,-15){$G'$}  
\put(-73,25){$1$} \put(-53,25){$2$} \put(-23,25){$3$} \put(-65,-10){$4$}
  \end{center}
  \caption{}
 \end{figure} 

\noindent 
White vertices may have extra edges which are not shown, and the invariant $\langle G\rangle$ is taken with respect to  a 3-valent vertex marked by a black square.

\begin{lem}  \label{lemnewsplit} Let $G_s$ be obtained from  $G$ by splitting the edge $3$  as shown above (middle). Then $\langle G_s\rangle_{1,2,3}= \Lef \langle G\rangle_{1,2,3}$.
\end{lem}
\begin{proof}
Let $f'_0,f'_1,f'_2,f'_3,f'_{123}$   denote the polynomials in $G_s$ with respect to the marked 3-valent vertex.  By \S\ref{sectlocal} (1) we have
$f'_0=f_0$, $f_1'=f_1$, $f_2'=f_2$, $f_3'=f_3+f_0\alpha_4$, $f_{123}'=f_{123}+ (f_1+f_2)\alpha_4$, where $f_0,f_1,f_2,f_3,f_{123}$ are the corresponding polynomials for $G$. It follows
immediately from the definitions that
$$\langle G_s\rangle_{1,2,3}=[f_0,\hspace{-0.5pt}f_1,\hspace{-0.5pt}f_2,\hspace{-0.5pt}f_3+f_0\alpha_4,\hspace{-0.5pt}f_{123}+(f_1+f_2)\alpha_4]
=\Lef [f_0,\hspace{-0.5pt}f_1,\hspace{-0.5pt}f_2,\hspace{-0.5pt}f_3,\hspace{-0.5pt}f_{123}]=\Lef\langle G\rangle_{1,2,3}.$$
\end{proof}

\begin{remark}  It can happen that two or more of the  edges $1,2,3$ of $G$ have two common endpoints (i.e., two or more of the  white vertices can coincide).
In this  degenerate case the invariant $\langle G \rangle$ is easily  expressed in terms of  graph polynomials $[\Psi_H]$, where $H$ is a  strict minor of $G$, by \S\ref{sectlocal} (2).
\end{remark}

\begin{lem} \label{lem3vtricase} Let $G'$ be as indicated above, with edges $1,2,3$ forming a three-valent vertex and edges $1,2,4$ forming a triangle, and let $G=G'\backslash \{4\}$.
Then 
$$\Lef \langle G'\rangle_{1,2,3} = (\Lef^2-\Lef) \langle G \rangle_{1,2,3} +[\Psi_{G,12} ] +[\Psi^3_{G,12}] - \Lef^{N_G-2} \ . $$
\end{lem}

\begin{proof} Since $G'$ has a three-valent vertex,  $\Psi_{G'}$ has the general shape given by lemma \ref{lemvtxstructure} with coefficients $f_0',f_1',f_2',f_3',f_{123}'$ where, 
by contraction-deletion:
$$f_0'= f_1+f_2+ f_0 \alpha_4\, ,  \quad  f_1'=f_1 \alpha_4 \, ,\quad   f_2'= f_2 \alpha_4  \, ,\quad f_3'=f_{123}+ f_3\alpha_4 \,  ,   \quad f_{123}' = f_{123} \alpha_4\, , $$
and $f_0,f_1,f_2,f_3,f_{123}$ are the corresponding structure constants for $G$.
On considering the two cases $\alpha_4=0$ and $\alpha_4\neq 0$ we find
$$[f_0',f_1',f_2',f_3',f_{123}' ] = (\Lef-1) [f_0,f_1,f_2,f_3,f_{123}] + [f_1+f_2,f_{123}]\ .$$
By definition \ref{defnfi}, $[f_1+f_2,f_{123}] = [\Psi^{3}_H, \Psi_{H,3}]$ where $H= G\q \{1,2\}$.   One concludes by applying  lemma \ref{lemlin} {\it (i)}.
\end{proof}

\begin{cor} \label{prophatreduction} Let $G'$ contain a split vertex as depicted in figure 3. Then
$$\Lef\langle G' \rangle_{1,2,4} =(\Lef^3-\Lef^2)\langle G\rangle_{1,2,3} + [\Psi_{G,2}]+ [\Psi^1_{G,2}] -\Lef^{N_G-2}\ .$$
\end{cor}
\begin{proof}  This follows from applying lemma $\ref{lem3vtricase}$ and then lemma $\ref{lemnewsplit}$ to figure 3.
\end{proof}

\subsection{Graphs of vertex width $\leq 3$} The notion of vertex width is  minor monotone.
\begin{lem} Let $G$ be a connected graph of vertex width $\leq n$, and let $H$ be any connected minor of $G$. Then the vertex width of $H$ is $\leq n$.
\end{lem}  
\begin{proof} Any ordering  on the edges of $G$ defines  a (strict) filtration  $G_i$ of subgraphs of $G$. This  induces a filtration  $H_i$ of subgraphs of  $H$ (which is not necessarily
strict any more). Clearly $|\hbox{vertices}(H_i \cap ( H\backslash H_i))| \leq |\hbox{vertices}(G_i \cap (G \backslash G_i)) |$.
\end{proof}

We give a constructive proof of the  following theorem (compare \cite{BrFeyn}, \S7.5) 
\begin{thm}  \label{thmvw3} If  $G$  has vertex-width at most 3,  then  $[\Psi_G]$ is a polynomial in $\Lef$.
\end{thm}

\begin{proof} A graph of vertex width $\leq 3$ comes with a filtration $G_i \subset G$ such that $G_i \cap (G\backslash G_i)$ has at most 3 vertices for all $i$.
For every minor $H$ of $G$ there is an induced ordering on its edges.
We say  that $H$ has an initial 3-valent vertex $v$, if $v$ is 3-valent and formed by the first three edges in $H$.

We show  that:
\begin{enumerate}
\item If  $G$ has an initial 3-valent vertex $v$, then $\langle G\rangle_{v}$ is  a  linear combination of  $[\Psi_H]$ and $\langle H\rangle_{v'}$ 
with coefficients in $\Z[\Lef]$, where $H$ are  strict minors of $G$, and $v'$ is an initial  3-valent vertex in $H$.
\item $[\Psi_G]$ is  a linear combination of $[\Psi_H]$ and $\langle H\rangle_{v'}$ with coefficients in $\Z[\Lef]$, where $H$ are strict minors of $G$, and $v'$ is an initial 3-valent
vertex  of $H$. 
\end{enumerate}
 These two facts, together with the fact that the vertex-width is minor monotone, are enough to prove the theorem. Note that if $H$ is not connected, then both
$[\Psi_H]$ and $\langle H\rangle_{v'}$ vanish.

First we show $(1)$. Consider  the subgraph  $G_4\subset G$ defined by its first four edges, and let $v$ be an initial  3-valent vertex in $G_4$. Suppose that $G_4\cap (G\backslash G_4)$
consists of exactly 3 distinct vertices (the degenerate cases where there are $\leq 2$  vertices are trivial by lemma \ref{SPlemma} and left to the reader).
Drawing the three vertices in white, we find ourselves in  the two cases denoted $G_s$ and $G'$ of figure 5 (up to renumbering of the edges).
Statement $(1)$ follows from lemmas \ref{lemnewsplit}  and \ref{lem3vtricase}.
 
Now we prove $(2)$. Consider the subgraph $G_5 \subset G$ defined by the first five edges of $G$. Assume that $G_5 \cap (G\backslash G_5)$ consists of three distinct
vertices $v_1,v_2,v_3$, since the degenerate cases where there are fewer than three vertices are again trivial by lemma \ref{SPlemma}.
If $G_5$ is not simple, then $[\Psi_G]$ trivially reduces to a linear combination of classes $[\Psi_H]$ where $H$  are strict minors of $G$, with coefficients in $\Z[\Lef]$.
  This follows from  lemma \ref{SPlemma}. If $G_5$ is  a simple graph,  then  up to renumbering of the edges, there are only two cases, shown  below:

 \begin{figure}[ht!]
  \begin{center}
    \epsfxsize=10.0cm \epsfbox{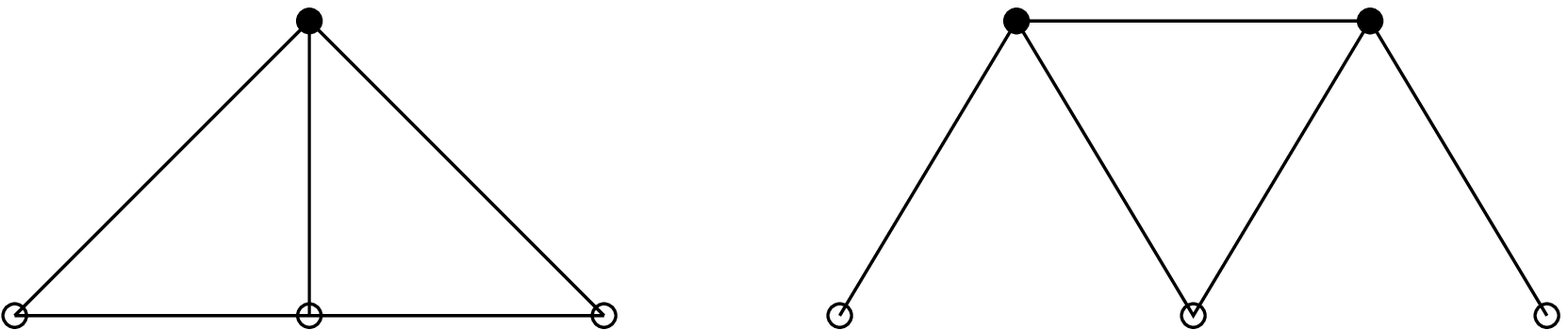}
 \put(-290,-10){$v_1$} \put(-235,-10){$v_2$} \put(-180,-10){$v_3$}
 \put(-137,-10){$v_1$} \put(-72,-10){$v_2$} \put(-10,-10){$v_3$}
  \label{figure3}
  \end{center}
\end{figure}

The left-hand figure is a split triangle and is covered by theorem \ref{thmsplittriangle}; the right-hand figure is a split vertex and is covered by theorem \ref{thmvertexsplit}.
Statement $(2)$ holds in both cases, which completes the proof.
\end{proof}

\subsection{Example 1: wheels with $n$ spokes} \label{sectwheels}
We use the previous results to compute the classes $[W_n]$ for all $n$, where $W_n$ denotes the wheel with $n$ spokes graph pictured below (left). Let $B_n$ denote the family
of graphs obtained by contracting a spoke of $W_n$, which have exactly $n$ vertices on the outer circle (right).

\begin{figure}[ht!]
 \begin{center}
    \leavevmode
    \epsfxsize=10.0cm \epsfbox{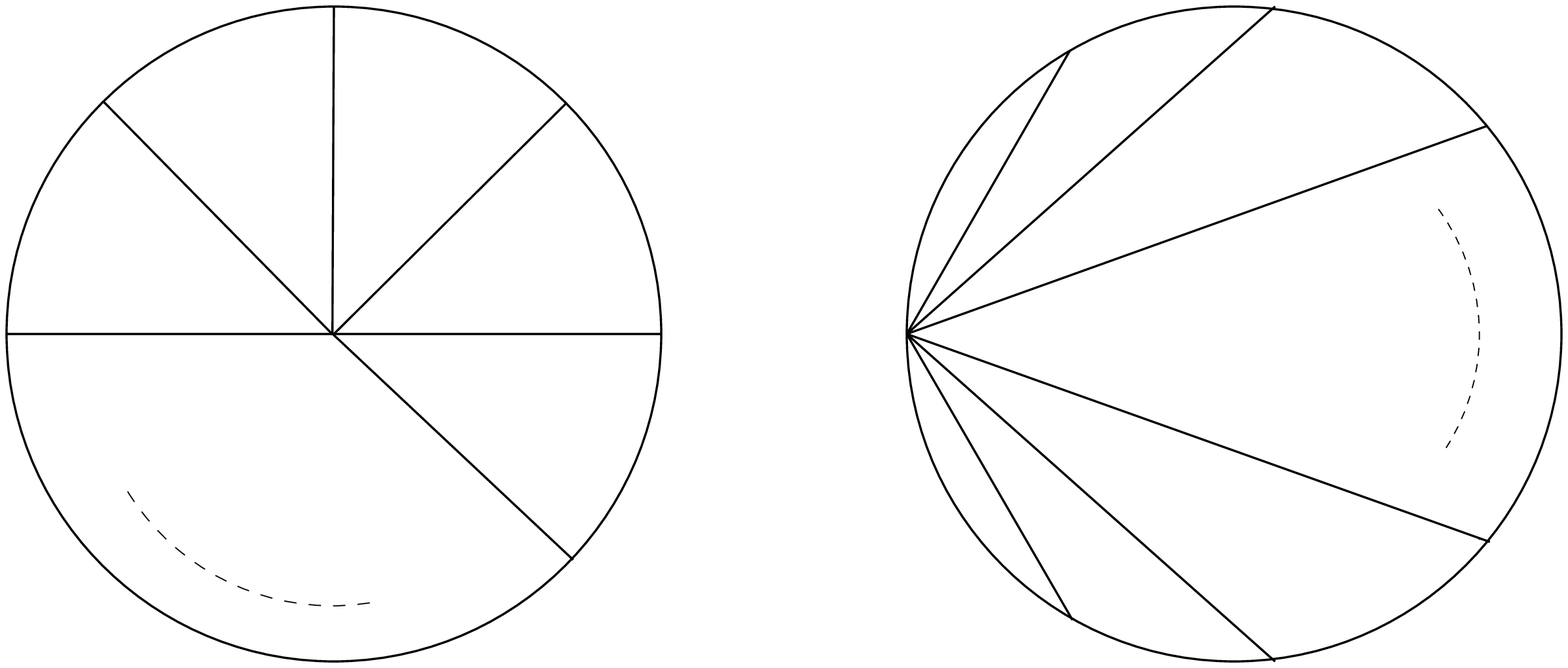}
  \put(-256,117){$4$}   \put(-208,120){$2$}  \put(-178,80){$5$}
  \put(-235,90){$1$}   \put(-200,75){$3$}  
  \put(-120,95){$e_1$}  \put(-97,95){$e_2$}
  \put(-282,105){$v_1$}   \put(-235,52){$v_2$}  \put(-163,60){$v_3$}
   \put(-234,125){$v_4$}   \put(-185,107){$v_5$} 
  \put(-315,10){\large{$W_n$}}  \put(-150,10){\large{$B_n$}}
  \end{center}
  \caption{The wheels with spokes graphs $W_n$, and a related family $B_n$ of series-parallel graphs.}
\end{figure}
The graphs $B_n$ are series-parallel reducible, so the classes $[B_n]$ can be computed using lemma \ref{SPlemma}.  This also follows from the results of \cite{AM}, theorem 5.10. 

\begin{lem} \label{Blem} Let us set $b_0=0$, $b_1=1$, and $b_n=[B_n]$ for $n\geq 2$.  If $B(t) = \sum_{n\geq 0 } b_nt^n$ is the generating series for the family of graphs $B_n$,
then we have
\begin{equation} \label{Bgenfun}
B(t) = \frac{(1 + \frac{\Lef t}{1-\Lef^2t})t}{1- (\Lef-1)(1+\Lef t)\Lef t }  \ .\end{equation}
\end{lem}
\begin{proof}
We refer to the two edges $e_1,e_2$ indicated on the diagram above. Since $e_1,e_2$ form a doubled edge, we have by $(\ref{Parallel})$:
$$[B_n]= (\Lef-2) [B_n \backslash e_1]+(\Lef-1)[B_n\backslash \{e_1,e_2\}] +\Lef [B_n \backslash e_1 \q e_2]+\Lef^{2n-3}$$
since $B_n$ has $2n-1$ edges. Now $B_n\backslash e_1$ is isomorphic to the graph obtained from $B_{n-1}$ by subdividing an outer edge, so $[B_n\backslash e_1 ]= \Lef[B_{n-1}]$ 
by $(\ref{Series})$. 
The graph $B_n\backslash \{e_1,e_2\}$ has an external leg, which
provides a factor of $\Lef$, leaving, as before,  a copy of $B_{n-2}$ with a subdivided outer edge. Thus $[B_n \backslash \{e_1,e_2\}]=\Lef^2 [B_{n-2}]$.
Finally, we have $B_n\backslash e_1 \q e_2 \cong B_{n-1}$, so we obtain
$$[B_n] = \Lef (\Lef-2) [B_{n-1}] +\Lef^2 (\Lef-1)[B_{n-2}] +\Lef [B_{n-1}]+\Lef^{2n-3}\ .$$
We deduce that for all $n\geq 4$ we have:
\begin{equation}\label{Brec} b_n = \Lef(\Lef-1) b_{n-1} + \Lef^2(\Lef-1) b_{n-2} +\Lef^{2n-3}\ .  \quad \quad \end{equation}
The constants $b_0,b_1$ are chosen such that the equation is valid for $n=2,3$, where $b_2=\Lef^2$  and $b_3 = \Lef^2 
(\Lef^2+\Lef-1)$ by direct computation. The formula for the  generating series then 
follows immediately from the recurrence relation $(\ref{Brec})$.
\end{proof}
One has $b_2=\Lef^2$, and  
{\small 
\begin{eqnarray}
b_3   =   \Lef^2 (\Lef^2+ \Lef-1 ) & , & b_4  =   \Lef^3 (\Lef^3 + 2\Lef^2-3\Lef+1)\ ,\nonumber \\
b_5  =  \Lef^5 (\Lef^3+ 3 \Lef^2-5\Lef+2) &  , &   b_6  =   \Lef^5 (\Lef^5 +4 \Lef^4 -7\Lef^3+2\Lef^2+2\Lef-1 ) \  .\nonumber
\end{eqnarray}
} 

Let $v_1,v_2,v_5$ denote any three vertices on $W_n$, joined by a three-valent vertex ($v_4)$ as shown in the diagram above. Let us write
$\langle W_n \rangle = \langle W_n \rangle_{v_4}$.
\begin{lem}  Let $\widehat{w}_n=0$ for $n\leq 2$ and  set  $\widehat{w}_n = \langle W_n\rangle$ for $n\geq 3$. Denote the corresponding ordinary generating series by
$\widehat{W}(t) = \sum_{n\geq 0 } \widehat{w}_n t^n$. Then
\begin{equation} \label{Whatgen}
\widehat{W}(t) = \frac{ (1+\Lef t)t \,   B(t) - \frac{ t^2}{1-\Lef^2 t}}{\Lef-(\Lef-1)\Lef^2t}  \ .\end{equation}
\end{lem}

\begin{proof} We use corollary \ref{prophatreduction}, applied to the graphs $G'= W_{n}$ with the edge and vertex labels on $G'$ as shown above. 
Then $G\cong  G' \backslash 1\q2 \cong W_{n-1}$. We have
$$\Lef \,\widehat{w}_n =(\Lef^3-\Lef^2) \, \widehat{w}_{n-1} + [\Psi^2_{G',45}]+  [\Psi^{12}_{G',45}] -\Lef^{2n-4} $$
since $W_n$ has $2n$ edges. Now $G'\backslash 2 \q \{4,5\}$ is isomorphic to $B_{n-1}$ and $G' \backslash \{1,2\} \q \{4,5\}$ gives the graph obtained from $B_{n-2}$ by subdividing
one outer edge. Therefore $ [\Psi^{12}_{G',45}]= \Lef [B_{n-2}]$ by lemma \ref{SPlemma}. We deduce that for all $n\geq 4$,
\begin{equation} \label{whatrec}
\Lef \,\widehat{w}_n =(\Lef^3-\Lef^2) \, \widehat{w}_{n-1} + b_{n-1}+  \Lef b_{n-2} -\Lef^{2n-4} \ .
\end{equation}
Using the fact that $\widehat{w}_3=1$ determines $\widehat{w}_n$ for $n=0,1,2$. The formula for the generating series $\widehat{W}(t)$ then follows 
immediately from $(\ref{whatrec})$.
\end{proof}

\begin{prop} Let $w_1=\Lef$, $w_2 = \Lef^3$, and  $w_n = [W_n]$  for $n\geq 3$. Let $W(t) =\sum_{n\geq 0 } w_n t^n$ be the generating function for the wheels with spokes graphs. Then 
\begin{equation} \label{Wgf}
W(t) =\frac{ (\Lef^4-\Lef^3)\,  \widehat{W}(t) +(\Lef-1) (1-\Lef^2t^2) B(t) +  \frac{ (1-\Lef t^2+\Lef^2t^2)}{(1-\Lef^2t)} }{1+(\Lef-1)t}\Lef t  
\end{equation}
where $B(t)$, $\widehat{W}(t)$ are defined above.
\end{prop}
\begin{proof}The graphs $W_1$ and $W_2$ are series-parallel and therefore $w_1$ and $w_2$ are given by lemma \ref{SPlemma}.
We apply theorem $\ref{thmvertexsplit}$ to the graph $G'= W_n$ with the labelling on its edges depicted above.  Since $G\cong W_{n-1}$, we deduce for all $n\geq 3$ that
$$[W_n] + (\Lef-\Lef^2)\big([\Psi_{G,2}] - [\Psi^{13}_{G,2}]\big)+(\Lef-1) [W_{n-1}]\hspace{-0.8pt}
=\hspace{-0.8pt} (\Lef^5-\Lef^4)\langle W_{n-1} \rangle   +\Lef^{2n-4}(\Lef^3+\Lef-1)$$
As before, $G\q \{2\} \cong B_{n-1}$, and $G\backslash \{1,3\}\q \{2\}$  is isomorphic to the graph obtained from $B_{n-3}$ by subdividing two  outer edges. It follows that
$[G\backslash \{1,3\}\q \{2\}]=\Lef^2 [B_{n-3}]$, giving
$$w_n + (\Lef-1) w_{n-1} + (\Lef-\Lef^2)\big(b_{n-1} - \Lef^2 b_{n-3}\big) = (\Lef^5-\Lef^4)\widehat{w}_{n-1}    +\Lef^{2n-4}(\Lef^3+\Lef-1)$$
The formula for the generating function follows from this.
\end{proof}

\begin{cor} Let $c_i(W_n)$ denote the coefficient of $\Lef^i$ in $[W_n]$.  Then 
$c_2(W_n) =-1$, $ c_{2n-1}(W_n) = 1$ and  $c_{2n-2}(W_n) =0$  for all $n$.  The outermost non-trivial coefficients are   $c_3(W_3) =1$, and 
$c_3(W_n) = n$ for all $n\geq 4$,  and
$ c_{2n-3} (W_n) = \binom{n}{2} $ for all $n\geq 3$.
\end{cor}

The following curious identity follows from the explicit description of $W(t)$:
$$[W_n] - [W_n\backslash O] - [W_n \q O] + [W_n\q I] = - \Lef^2 (\Lef-1)^{n-2}\ ,$$
where $O$ denotes any outer edge of $W_n$ (on the rim of the wheel), and $I$ denotes any internal edge or spoke. The combinatorial reason for this is not clear.

\begin{remark} The polynomials $w_n$ should have  equivariant versions with respect to the symmetry group of $W_n$.
 Computing these explicitly  would be relevant to  computing  the full cohomology of the graph hypersurface complement of $W_n$ and, what one ultimately wants: the action of
 the motivic Galois group.
\end{remark}
\noindent The first few values of the polynomials $w_n$ are as follows:
{\small
\begin{eqnarray}
w_3 & = & \Lef^2 (\Lef^3+\Lef-1) \nonumber \\ 
w_4 & = &  \Lef^2 (\Lef^5   +3\Lef^3-6\Lef^2 +4\Lef-1) \nonumber        \\  
w_5  & = &  \Lef^2 (\Lef^7  +6\Lef^5-15\Lef^4 +16\Lef^3  -11 \Lef^2 +5\Lef-1) \nonumber   \\
w_6 & = & \Lef^2(\Lef^9 +10 \Lef^7 -29 \Lef^6+37\Lef^5  -33\Lef^4+26\Lef^3-16\Lef^2+6\Lef-1) \nonumber          \\
w_7 & = & \Lef^2 (\Lef^{11}+15 \Lef^9-49 \Lef^8+71 \Lef^7-70 \Lef^6+64 \Lef^5-57 \Lef^4+42 \Lef^3-22 \Lef^2+7 \Lef-1) \nonumber 
\end{eqnarray}}
Note that the wheels $W_n$ are the unique
infinite family of graphs whose residue
can be calculated, namely: 
$I_{W_n} = \binom{2n-1}{n-1} \zeta(2n-3)$ for $n\geq 3$. 
One of the main results of \cite{B-E-K} is that
\begin{equation}\label{cohomBEK}  H^{2n-1}_c (\Pro^{2n-1}\backslash X_{W_n}) \cong \Q(-2)\end{equation}
and that $H^{2n-1} (\Pro^{2n-1}\backslash X_{W_n})$ is generated by the class of the integrand of  $I_G$. 
It would be interesting to relate their proof to the above computation which   gives    $c_2(W_n) =-1$. 

\subsection{Example 2:  Zig-zag graphs} \label{sectZigs}
The second application of the previous results  is  to compute the classes $[Z_n]$ for all $n$, where $Z_n$ denotes the
family of zig-zag graphs with $n$ loops pictured below (left).
 Let $\overline{Z}_n$ denote the family
of graphs obtained by doubling the   edge  `2' as shown on the right. Note that $Z_3=W_3$.

\begin{figure}[ht!]
 \begin{center}
    \leavevmode
    \epsfxsize=11.8cm \epsfbox{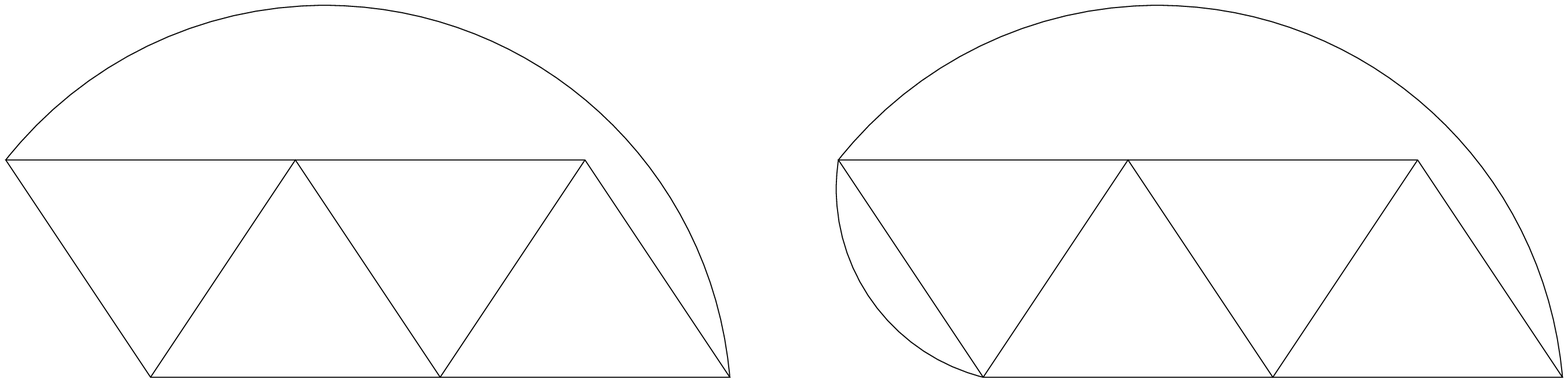}
 \label{ZZ}
   \put(-340,75){$Z_5$} 
  \put(-320,20){$2$} \put(-290,20){$3$}
  \put(-308,50){$1$}   \put(-260,70){$4$}  \put(-280,-10){$5$}
    \put(-347,41){$v_4$}  \put(-315,-6){$v_5$} \put(-250,-6){$v_3$}\put(-278,50){$v_2$}
    \put(-170,75){$\overline{Z}_5$}\put(-183,-6){$v_1$}
  \end{center}
  \caption{The zig-zag graphs $Z_n$.}
\end{figure}
The graphs $Z_n$ are primitive-divergent graphs in $\phi^4$ theory for all $n\geq 3$.
Let $z_0=0$,  $z_1=\Lef+1$, $z_2= \Lef^3$, and $z_n =[Z_n]$ for all $n\geq 3$. Likewise, set 
$\overline{z}_0=1$, $\overline{z}_1 = \Lef^2$, $\overline{z}_2 = \Lef^4+\Lef^3-\Lef^2$, and $\overline{z}_n = [\overline{Z}_n]$ for all $n\geq 3$.
Denote the corresponding generating series by $Z(t)$ and $\overline{Z}(t)$.
\begin{lem} A straightforward application of the series-parallel operations gives
\begin{equation} \label{ZR1}
\overline{z}_n = (\Lef-2) z_n + (\Lef-1)\Lef^2 \overline{z}_{n-2} + \Lef \overline{z}_{n-1} + \Lef^{2n-1} \quad n\geq 1
\end{equation}
\end{lem}
\begin{proof} If $e_1,e_2$ denote the two doubled edges, then  this follows from  the parallel reduction $(\ref{Parallel})$ on noting that
$\overline{Z}_n\backslash e_1 \cong Z_n$, $\overline{Z}_n \backslash e_1 \q e_2 \cong \overline{Z}_{n-1}$, and
that $\overline{Z}_n \backslash \{e_1,e_2\}$ is isomorphic to the graph obtained from  $Z_{n-2}$ by subdividing two  edges, whose class is $\Lef^2 z_{n-2}$ by two applications of 
$(\ref{Series})$. 
\end{proof}

We next want to compute recursion relations for the numbers $z_n$ by  considering the split vertex shown above in the figure (left).
Let us set $\widehat{z}_n = 0 $ for $ n <3$,   $\widehat{z}_n =  \langle Z_n \rangle_{v_4}  $ for all $n\geq 3$, and let $\widehat{Z}(t)$ be the corresponding generating series.
Let $ZB_n$ denote the family of graphs depicted below with $n$ vertices. A trivial argument along the lines of lemma $\ref{Blem}$ shows that $[ZB_n]=b_n$, with
generating series $B$.

\begin{figure}[ht!]
 \begin{center}
    \leavevmode
    \epsfxsize=5cm \epsfbox{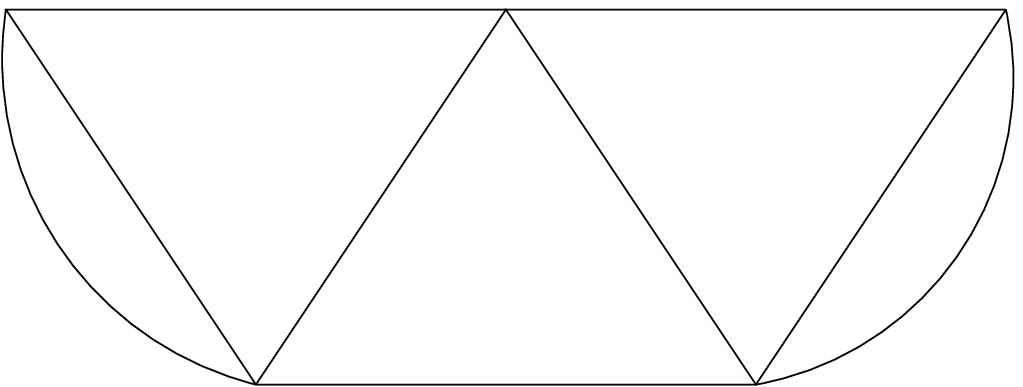}
   \put(-190,40){$ZB_5$} 
  \end{center}
\end{figure}

\begin{lem} The recurrence relation given in  theorem $\ref{thmvertexsplit}$ translates as:
\begin{equation}\label{ZR2}
z_n+(\Lef-\Lef^2)(\overline{z}_{n-2}-\Lef^2 b_{n-3}) +(\Lef-1)z_{n-1} = (\Lef^5-\Lef^4)\widehat{z}_{n-1}+\Lef^{2n-4}(\Lef^3+\Lef-1)
\end{equation}
for $n\geq 2$.  The recurrence relation of  corollary $\ref{prophatreduction}$ yields the relation
\begin{equation} \label{ZR3}
\Lef \widehat{z}_n= (\Lef^3-\Lef^2)\widehat{z}_{n-1} +\overline{z}_{n-2} +\Lef b_{n-2}-\Lef^{2n-4} \ ,  \qquad n\geq 2
\end{equation}
\end{lem}
\begin{proof} Let $G'= Z_n$, and apply theorem $\ref{thmvertexsplit}$ to $G'$ with the edge numbering shown above. Then $G\cong Z_{n-1}$, 
$G'\backslash 2 \q \{4,5\} \cong \overline{Z}_{n-2}$, and $G\backslash \{1,3\} \q 2$ is isomorphic to the graph obtained from $ZB_{n-3}$ by subdividing two edges.
It follows from  $(\ref{Series})$ that $[\Psi_{G\backslash \{1,3\} \q 2}] = \Lef^2 b_{n-3}$, which yields the first equation.
The second equation follows from corollary $\ref{prophatreduction}$, since  $G'\backslash \{1,2\}\q \{4,5\} $ is isomorphic to the graph obtained from $ZB_{n-2}$
by subdividing one edge, whose polynomial is $\Lef b_{n-2}$ by $(\ref{Series})$.
\end{proof}

Equations $(\ref{ZR1})$, $(\ref{ZR2})$, $(\ref{ZR3})$ imply the following identities of generating series:
\begin{eqnarray*}
&&\hspace*{-1cm}{[}1-\Lef t + \Lef^2(1 \!- \!\Lef)t^2]\, \overline{Z} -(\Lef-2) Z - \Lef \, R  = 1+(2\!-\!\Lef)t \\
&&\hspace*{-1cm}{[} \Lef^3-(\Lef ^5-\Lef ^4) t ]\, \widehat{Z}  - \Lef ^2t^2\big( \overline{Z}+\Lef  B \big)  + R =t \\
&&\hspace*{-1cm}{[} 1+(\Lef \!-\!1)t] Z-(\Lef ^5\!-\!\Lef ^4)t \widehat{Z} + (\Lef \!-\!\Lef ^2) t^2 ( \overline{Z}- \Lef ^2 t B)-(\Lef ^3\!+\!\Lef -1)tR =(\Lef \!+\!1) t
\end{eqnarray*}
 in three unknowns,  $Z$, $\overline{Z}$, and $\widehat{Z}$, where
  $R=R(t) = \frac{t}{(1-\Lef^2t)}.$  These equations are easily solved using  the expression for $B$ $(\ref{Bgenfun})$.
   In particular, we  obtain an explicit  formula for the  generating series  for the zig-zag graphs:
  \begin{equation} \label{zzgen}  {1\over 1-\Lef^2 t}-{   \M^2\Lef^2t^3P(\Lef , \Lef t) \over \big( 1-    \M \Lef t - \M \Lef^2t^2 \big)^2\big(1-t-\M^2\Lef t^2\big)}
 \end{equation}
where $\M= \Lef -1 $ and the polynomial $P$ is defined by:
\begin{eqnarray} P(x,y) &=& 
x(x-1)^3 y^5+(x-1)(2x^3-3x^2+x+1)y^4 \nonumber \\
& &   +\quad (x-1)(x^3-3x^2+2x+1)y^3-(3x^3-3x^2+2)y^2 \nonumber \\
&& -\quad (x^3-x^2+1)y+x^2+x+1  \nonumber \end{eqnarray}  
The coefficient of $t^n$ in $(\ref{zzgen})$ is $z_n$ for $n\geq 3$ (see below).
This is to our knowledge the only explicit formula for the class in the Grothendieck ring of a family of primitive-divergent graphs in $\phi^4$.
 From this formula   one  obtains:

\begin{cor} Let $c_i(Z_n)$ denote the coefficient of $\Lef^i$ in $[Z_n]$.  Then 
$c_2(Z_n) =-1$, $ c_{2n-1}(Z_n) = 1$ and  $c_{2n-2}(Z_n) =0$  for all $n$.  The outermost non-trivial terms are   $c_3(Z_3) =1$, and 
$c_3(Z_n) = 8-n$ for all $n\geq 4$,  and
$ c_{2n-3} (W_n) = 2n-5 $ for  $n\geq 3$.
\end{cor}
In the case of the zig-zags, the analogous result to $(\ref{cohomBEK})$ was proved 
for $n\geq5$ by Doryn in his thesis \cite{Doryn}. It states that 
$\gr^W_{min} H^{2n-1}_c (\Pro^{2n-1}\backslash X_{Z_n}) \cong \Q(-2), $
which should again be related to our proof that $c_2(Z_n)=-1$.

For small $n$, we have:
{\small
\begin{eqnarray}
z_3 & = & \Lef ^2 (\Lef ^3+\Lef -1) \nonumber \\ 
z_4 & = &    \Lef ^2( \Lef ^5+  3\Lef ^3 -6\Lef ^2 +4\Lef -1)  \nonumber        \\  
z_5  & = &  \Lef ^2  (\Lef ^7+5 \Lef ^5-10 \Lef ^4+7 \Lef ^3-4 \Lef ^2+3 \Lef -1) \nonumber   \\
z_6 & = & \Lef ^2   (\Lef ^9+7 \Lef ^7-12 \Lef ^6-2 \Lef ^5+16 \Lef ^4-12 \Lef ^3+2 \Lef ^2+2 \Lef -1)  \nonumber    \\ 
z_7& = & \Lef ^2 (\Lef ^{11}+9\Lef ^9-13\Lef ^8-18\Lef ^7+55\Lef ^6-58\Lef ^5+41\Lef ^4-23\Lef ^3+7\Lef ^2+\Lef -1)  \nonumber       
\end{eqnarray}
}
Note that explicit results for the zig-zag periods were conjectured in \cite{BK}. Remarkably they are a rational
multiple of $\zeta(2n-3)$.

\newpage
\section{Non-Tate counter-examples at 8  loops} 
We use the denominator-reduction method to derive some non-Tate counter-examples to Kontsevich's conjecture at 8 and 9 loops.
\subsection{Combinatorial reductions} \label{sectCombinRed} In order to compute the $c_2$-invariant of an 8-loop graph $G$, we proceed  in two simpler steps.
The following lemmas will be applied to  the main counterexample, depicted in figure 8 below.

Suppose that  $G$ is any connected graph with the shape depicted below, where the white vertices $A,B,C,D$ may have anything attached to 
them. Let $H$ be the minor obtained from $G$ by deleting the edges $2$ and $4$, and contracting  $3$ and $6$.
\begin{center}
\fcolorbox{white}{white}{
  \begin{picture}(224,132) (37,-75)
    \SetWidth{1.0}
    \SetColor{Black}
    \Vertex(105,44){3}
    \SetWidth{2.0}
    \Line(105,44)(201,44)
    \Line(105,44)(137,-8)
    \Line(201,44)(169,-8)
    \Line(201,44)(233,-8)
    \Line(105,44)(73,-8)
    \Line(73,-8)(137,-8)
    \SetWidth{1.0}
    \Vertex(201,44){3}
    \SetWidth{2.0}
    \Arc(137,-8)(3,270,630)
    \Arc(169,-8)(3,270,630)
    \Arc(233,-8)(3,270,630)
    \Arc(73,-8)(3,270,630)
    \Text(67,-25)[lb]{{\Black{$A$}}}
    \Text(131,-25)[lb]{{\Black{$B$}}}
    \Text(169,-25)[lb]{{\Black{$C$}}}
    \Text(226,-25)[lb]{{\Black{$D$}}}
    \Text(35,30)[lb]{\Large{\Black{$G$}}}
    \Text(102,-20)[lb]{{\Black{$1$}}}
    \Text(124,19)[lb]{{\Black{$2$}}}
    \Text(79,19)[lb]{{\Black{$3$}}}
    \Text(145,51)[lb]{{\Black{$4$}}}
    \Text(175,19)[lb]{{\Black{$5$}}}
    \Text(230,19)[lb]{{\Black{$6$}}}
    \Arc(73,-71)(3,270,630)
    \Arc(73,-71)(3,270,630)
    \Arc(73,-71)(3,270,630)
    \Arc(137,-71)(3,270,630)
    \Arc(233,-71)(3,270,630)
    \Arc(169,-71)(3,270,630)
    \Line(73,-71)(137,-71)
    \Line(169,-71)(233,-71)
    \Text(226,-60)[lb]{{\Black{$D$}}}
    \Text(169,-60)[lb]{{\Black{$C$}}}
    \Text(131,-60)[lb]{{\Black{$B$}}}
    \Text(67,-60)[lb]{{\Black{$A$}}}
     \Text(102,-65)[lb]{{\Black{$1$}}}
 \Text(200,-65)[lb]{{\Black{$5$}}}
    \Text(34,-64)[lb]{\Large{\Black{$H$}}}
  \end{picture}
}
\end{center}
\begin{lem} \label{lemfirstred}
 Let $G, H$ be as above. Then
 ${D}^6_G(1,2,3,4,5,6) =\pm \Psi_{H}^{1,5} \Psi^5_{H,1}.$
\end{lem}

\begin{proof}  The proof is by direct computation of resultants, using the identities between Dodgson polynomials which follow from the existence of local stars and triangles.
Since the edges $\{1,2,3\}$ form a triangle,  we know from  example \ref{example5split} that
$${}^5 \Psi_G(1,2,3,4,5) = \pm\Psi^{123,345}_G \Psi^{4,5}_{G\backslash 2  \q 
\{1,3\}}\ .$$
Since $\{2,3,4\}$ forms a three-valent vertex,  we have $\Psi^{123,345}= \Psi^{1,5}_{G\backslash \{2,4\}\q 3}$ by the last equation in \S\ref
{sectlocal}, (3).
By contraction-deletion,  this last term is also $\Psi^{14,45}_{G\backslash 2\q 3}$, giving
$${}^5 \Psi_G(1,2,3,4,5) = \pm\Psi^{14,45}_{G_2}  \Psi^{4,5}_{G_2,1}\ ,$$
where $G_2$ is the minor  $G\backslash 2 \q 3$ with the induced numbering of its edges. 
Now take the resultant with respect to edge $6$. Since $\{4,5,6\}$ forms a three-valent vertex in $G_2$, it follows that
$\Psi_{G_2}^{146,456}=0 $ by the vanishing property for vertices. Thus we have
 $$[\Psi^{14,45}_{G_2},  \Psi^{4,5}_{G_2,1}]_{6}= \pm\Psi^{14,45}_{G_2,6} \Psi^{46,56}_{G_2,1}\ .$$
 Again, since $\{4,5,6\}$ is a three-valent vertex, $\Psi^{46,56}_{G_2,1} =  \Psi^{46,56}_{G_2\q 1} =  \Psi^{45}_{G_2\q
\{1, 6\}}= \Psi^{5}_{G_2\backslash 4\q\{1, 6\}}$, 
 where the first and third equality are  contraction-deletion relations. We have:
 $$[\Psi^{14,45}_{G_2},  \Psi^{4,5}_{G_2,1}]_{6}=\pm\Psi^{1,5}_{G_2\backslash 4 \q 6} \Psi^{5}_{G_2\backslash 4\q 6,1}$$ 
The left-hand side is  equal to $\pm {D}^6_G(1,2,3,4,5,6)$ by definition, and the minor $G_2\backslash 4\q 6$ is exactly $H$, which completes the proof.
\end{proof}

\begin{figure}[ht!]
 \begin{center}
    \leavevmode
    \epsfxsize=5cm \epsfbox{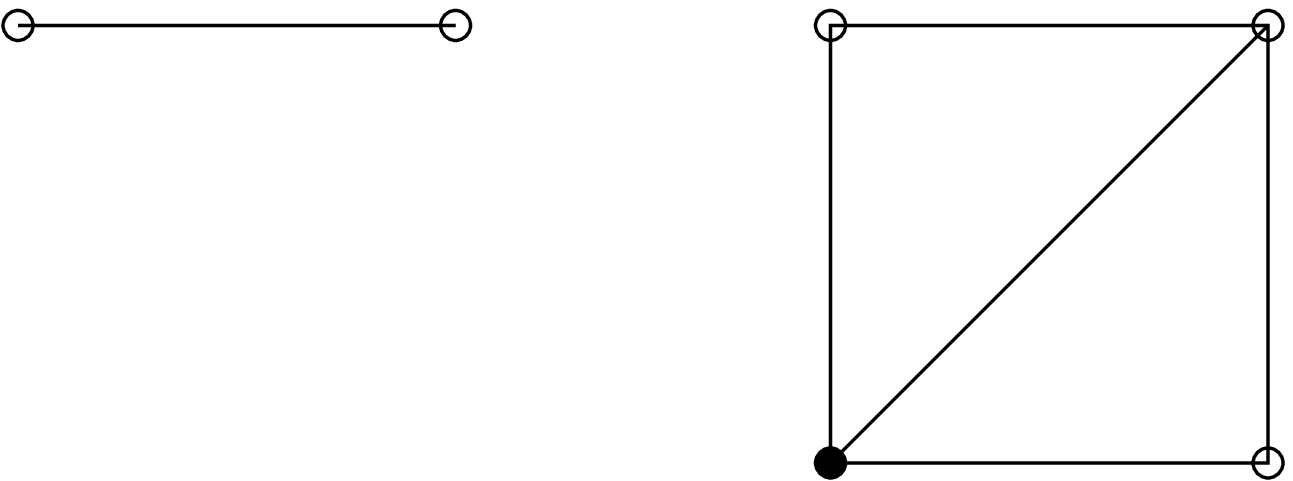}
  \put(-123,42){
$1$}\put(-35,42){
$5$}
  \put(-57,23){$8$} \put(-35,23){$7$} \put(-1,23){$10$}
  \put(-30,5){$9$}
    \put(-100,0){\large{$H$}}
  \end{center}
 \end{figure}

\begin{lem} \label{lemsecondred}  Now let $H$ be a graph with the general shape depicted above. The  denominator reduction, applied five times to 
$\Psi_H^{1,5} \Psi^5_{H,1}$ with respect to the edges $7,8,9,10$
is $\pm\Psi^{15,78}_A \Psi_B$, where $A= H\backslash \{10\}\q 9$ and $B=H\backslash \{5,7,9\}\q \{1,8,10\}$.
\end{lem}

\begin{proof}
By the second Dodgson identity, $[\Psi_H^{1,5},\Psi^5_{H,1}]_7 = \pm \Psi_H^{57,15}\Psi^{5,7}_{H,1}$. Applying the first Dodgson identity, we then  get 
 $[ \Psi_H^{57,15},\Psi^{5,7}_{H,1}]_8= \pm \Psi_H^{15,78}\Psi^{58,57}_{H,1}$. Now,
 $$[\Psi_H^{15,78},\Psi^{58,57}_{H,1}]_9 = -\Psi^{15,78}_{H,9} \Psi^{579,589}_{H,1}, $$
 by definition of the resultant, using the fact that $\Psi^{159,789}_{H}=0$, by the vanishing property for vertices applied to the three-valent vertex 7,8,9.
 Once more, by the vanishing property applied to the triangle $7,9,10$, we have
$\Psi^{15,78}_{H,9X}=0$ where $X$ denotes the edge 10, and therefore
$$[\Psi^{15,78}_{H,9} ,\Psi^{579,589}_{H,1}]_{10} =  \Psi^{15X,78X}_{H,9} \Psi^{579,589}_{H,1X}\ .$$
By contraction-deletion, the first factor is $\Psi^{15,78}_A$, and the second is
$\Psi^{7,8}_{H'}$ where $H'=H \backslash \{5,9\}\q \{1,10\}$. In this latter graph, $7,8$ form a 2-valent vertex, and so
$\Psi_{H'}^{7,8}=\Psi^7_{H',8} = \Psi_{H'\backslash 7\q 8}=\Psi_B$. 
\end{proof}

\subsection{An eight-loop counter-example} \label{sect8loopcounterNP}
Let $G_8$ be the eight-loop primitive-di\-ver\-gent $\phi^4$  graph with vertices numbered $1,\ldots, 9$ and (ordered) edges $e_1,\ldots, e_{16}$ defined by
\begin{equation} \label{edges8loopcounterNP}
34,14,13,12,27,25,58,78,89,59,49,47,35,36,67,69\ , 
\end{equation}
where $ij$ denotes an edge connecting vertices $i$ and $j$. 
\begin{figure}[ht!]
 \begin{center}
    \leavevmode 
    \epsfxsize=5cm \epsfbox{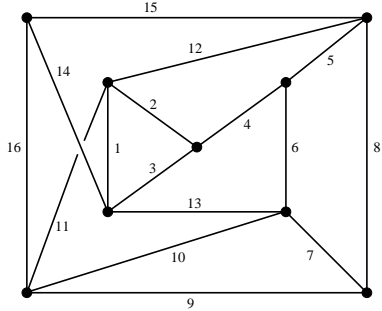}
    \end{center}
    \caption{The graph $G_8$}
      \end{figure}
      
This graph  is isomorphic to $P_{8,37}$ minus vertex 3 or 5 in the census \cite{SchnetzCensus}. It has 3785 spanning trees. The first six edges form precisely 
the configuration depicted in lemma $\ref{lemfirstred}$, and  we can subsequently apply lemma $\ref{lemsecondred}$ to reduce the  next four edges.
 A further reduction  with respect to edge $11$ gives the following corollary.

\begin{cor}
Let $G_8$ be the 8-loop graph defined by $(\ref{edges8loopcounterNP})$. Then 
$$D^{11}_{G_8} (e_1,\ldots, e_{11}) = \det 
\left(
\begin{array}{cc}
\Psi^{15,78}_{A\backslash 11} &  \Psi_{B\backslash 11}    \\
 \Psi^{15,78}_{A\q 11} &  \Psi_{B \q 11}    
\end{array}
\right)\ ,
 $$
where $A,B$ are depicted below.
\end{cor}

\vspace{-0.3in}

\begin{figure}[ht!]
 \begin{center}
    \leavevmode
    \epsfxsize=9.0cm \epsfbox{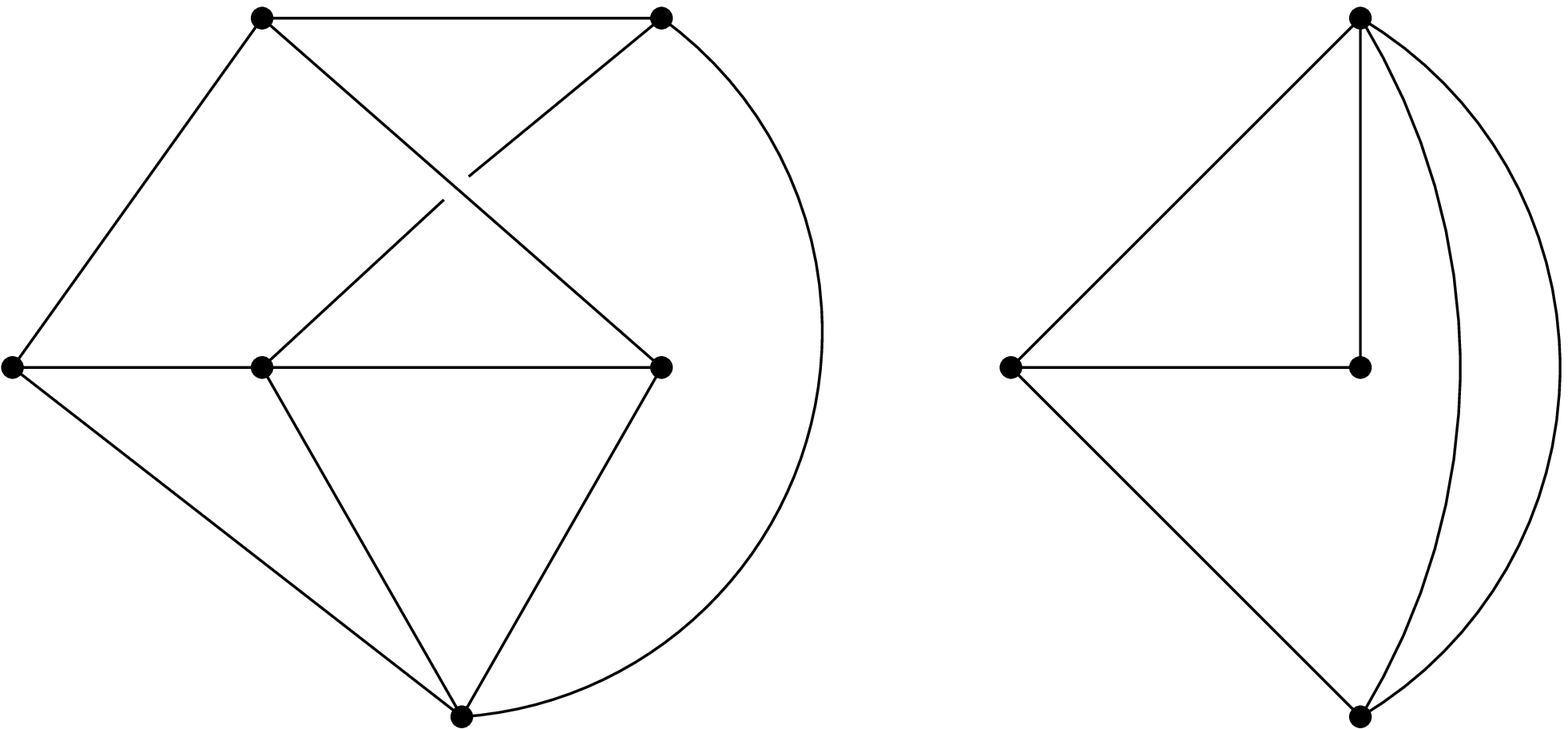}
 \put(-245,97){$14$}  \put(-212,95){$13$} 
 \put(-212,75){$12$} 
 \put(-245,65){$15$} 
 \put(-245,30){$16$} \put(-198,35){$8$}\put(-170,35){$7$}\put(-185,48){
$5$}\put(-185,120){
$1$}
  \put(-137,60){$11$}  \put(-75,20){$16$}\put(-20,50){$13$}\put(0,50){$11$} \put(-73,65){$15$} \put(-70,100){$14$}\put(-50,80){$12$}
  \put(-275,10){\large{$A$}}  \put(-105,10){\large{$B$}}
  \end{center}
\end{figure}

\vspace{-0.1in}
The polynomial $D^{11}_{G_8}(e_1,\ldots, e_{11})$ is irreducible, so to proceed further in the reduction, observe that $A$ and $B$ 
 have a common minor $\gamma=B\backslash\{11\}\q \{12,13\}$ which is the sunset graph on 2 vertices and 3 edges $14,15,16$. Its graph polynomial is
 $$\Psi_{\gamma}=  \alpha_{14}\alpha_{15}+\alpha_{15}\alpha_{16}+\alpha_{14}\alpha_{16}\ .$$
By direct computation, one verifies that
\begin{eqnarray} \label{detequations}
\Psi^{15,78}_{A\backslash 11} & = & -\alpha_{13}\alpha_{15}   \\
  \Psi^{15,78}_{A\q 11}  & = & \alpha_{12}(\Psi_{\gamma}+\alpha_{13}\alpha_{16}) \nonumber \\
 \Psi_{B\backslash 11}& = &  \Psi_{\gamma}+\alpha_{12}\alpha_{13}+\alpha_{16}\alpha_{12}+\alpha_{14}\alpha_{12}+\alpha_{15}\alpha_{13}+\alpha_{14}\alpha_{13} \nonumber \\
\Psi_{B \q 11}    & = &  \alpha_{13}( \Psi_{\gamma}+\alpha_{16} \alpha_{12}+ \alpha_{14} \alpha_{12})\ .  \nonumber 
\end{eqnarray}
By theorem \ref{corDenom}, $c_2(G_8)_q\equiv -[D^{11}_{G_8}(e_1,\ldots, e_{11})]_q\mod q$.
We can eliminate a further variable by exploiting the homogeneity of $D^{11}_{G_8}$ (or $\Psi_{G_8}$). The affine complement of the zero locus of a homogenous polynomial $F$
admits a $\mathbb{G}_m$ action    by scalar diagonal  multiplication of the coordinates. For any coordinate $\alpha_e$, we therefore  have
$$[F]_q= [F, \alpha_e]_q + (q-1) [F, \alpha_e-1]_q$$
\begin{lem} $[D^{11}_{G_8},\alpha_{16}]_q$ is a polynomial in $q$.
\end{lem}
\begin{proof} By inspection of $(\ref{detequations})$, setting $\alpha_{16}=0$ in the definition of $D^{11}_{G_8}$ causes the terms $\alpha_{14}\alpha_{15}$ to factor out.
The other factor is  of degree at most one in $\alpha_{14}$ and $\alpha_{15}$, and by a simple application of lemma \ref{lemlin} is therefore a polynomial in $q$.  
\end{proof}

We will  henceforth work on the hyperplane  $\alpha_{16}=1$. 
Now we may  scale $\alpha_{12}$ and $\alpha_{13}$ by $\Psi_{\gamma}$,  which has the effect of replacing $D^{11}_{G_8}$ with $\widetilde{D}$ given by formally setting
$\Psi_{\gamma}$ to be 1 and $\alpha_{12}\alpha_{13}$ to be
$\alpha_{12}\alpha_{13}\Psi_\gamma$ in the previous equations. Since this transformation is  an isomorphism on the complement of $V(\Psi_{\gamma})$, we have
\begin{equation}\label{9}
[D^{11}_{G_8} ]_q  -[\Psi_\gamma,D^{11}_{G_8}]_q  = [\widetilde{D} ]_q - [\Psi_{\gamma}, \widetilde{D}]_q 
 \ . 
\end{equation}
\begin{lem} $[\widetilde{D}]_q $ and $[\Psi_\gamma,D^{11}_{G_8}]_q$ are constant modulo $q$.
\end{lem}
\begin{proof} By inspection of $(\ref{detequations}),$ it is clear that the determinant $\widetilde{D}$ is of degree one in the variables $\alpha_{14}$ and $\alpha_{15}$.
Applying lemma \ref{lemlin} {\it (i)} twice,
it follows that the class of $\widetilde{D}$ modulo $q$ is equal to the class modulo $q$ of its coefficient of  $\alpha_{14}\alpha_{15}$, and this is 
 $\alpha_{12}\alpha_{13}( \alpha_{13}\alpha_{12}+\alpha_{13}+\alpha_{12})$, which gives a polynomial in $q$.
Likewise, a straightforward calculation using  $(\ref{detequations})$ shows that the intersection $V(\Psi_{\gamma}, D^{11}_{G_8})$ is union of intersections of hypersurfaces
of degree at most 2 and linear in every variable, which can be treated using lemma $\ref{lemlin}$ with components of small degree.
\end{proof}

It remains to compute $[\Psi_{\gamma},\widetilde{D}]_q$, which is given mod $q$ by  the resultant $[\Psi_{\gamma},\widetilde{D}]_{14}$.
Explicitly, it is  the polynomial:
\begin{equation}\label{penultimate}
\alpha_{12}+\alpha_{12}\alpha_{15}+\alpha_{13}\alpha_{12}^2+\alpha_{12}^2+\alpha_{13}\alpha_{12}+\alpha_{15}\alpha_{13}\alpha_{12} \end{equation}
$$ + \alpha_{13}^2\alpha_{15}+\alpha_{13}^2\alpha_{15}^2+\alpha_{13}^2\alpha_{15}\alpha_{12}+\alpha_{13}^2\alpha_{15}^2\alpha_{12}+\alpha_{15}^2\alpha_{13}\alpha_{12}\ .$$
A final innocuous change of variables  $\alpha_{13}\mapsto \alpha_{13}/ (\alpha_{15}+1)$
can be handled as in the previous case (\ref{9}) and reduces this equation to degree 4. Setting $a=\alpha_{13}+1$, $b=\alpha_{12}+1$, $c=\alpha_{15}$
leads to the equation
$$J = a^2bc-ab-ac^2-ac+b^2c+ab^2+abc^2-abc$$
which defines a singular surface in $\mathbb{A}^3$.  In conclusion
\begin{equation} \label{c2ofG8} c_2(G_8) \equiv c-[J]_q \mod q
\end{equation} for some constant $c \in 
\Z$. Chasing the  constant terms in the above gives  $c=2$.
Note that
$G$ has vertex-width 4 (realised by  a different ordering on the edges from the one given above). 
See also  \cite{SchnetzFq} for the complete computer-reduction of a graph in the same completion class as this one.
The proof that this is  a counter-example continues in \S\ref{sect7}, where we study the point counting function of  $V(J)$ in detail.

\subsection{A planar counter-example}\label{sectplanarcounter} Consider the  planar graph $G_9$  with nine loops and eighteen edges below. It is primitive-divergent and in $\phi^4$ theory.

\begin{figure}[ht!]
 \begin{center}
    \leavevmode 
    \epsfxsize=6.0cm \epsfbox{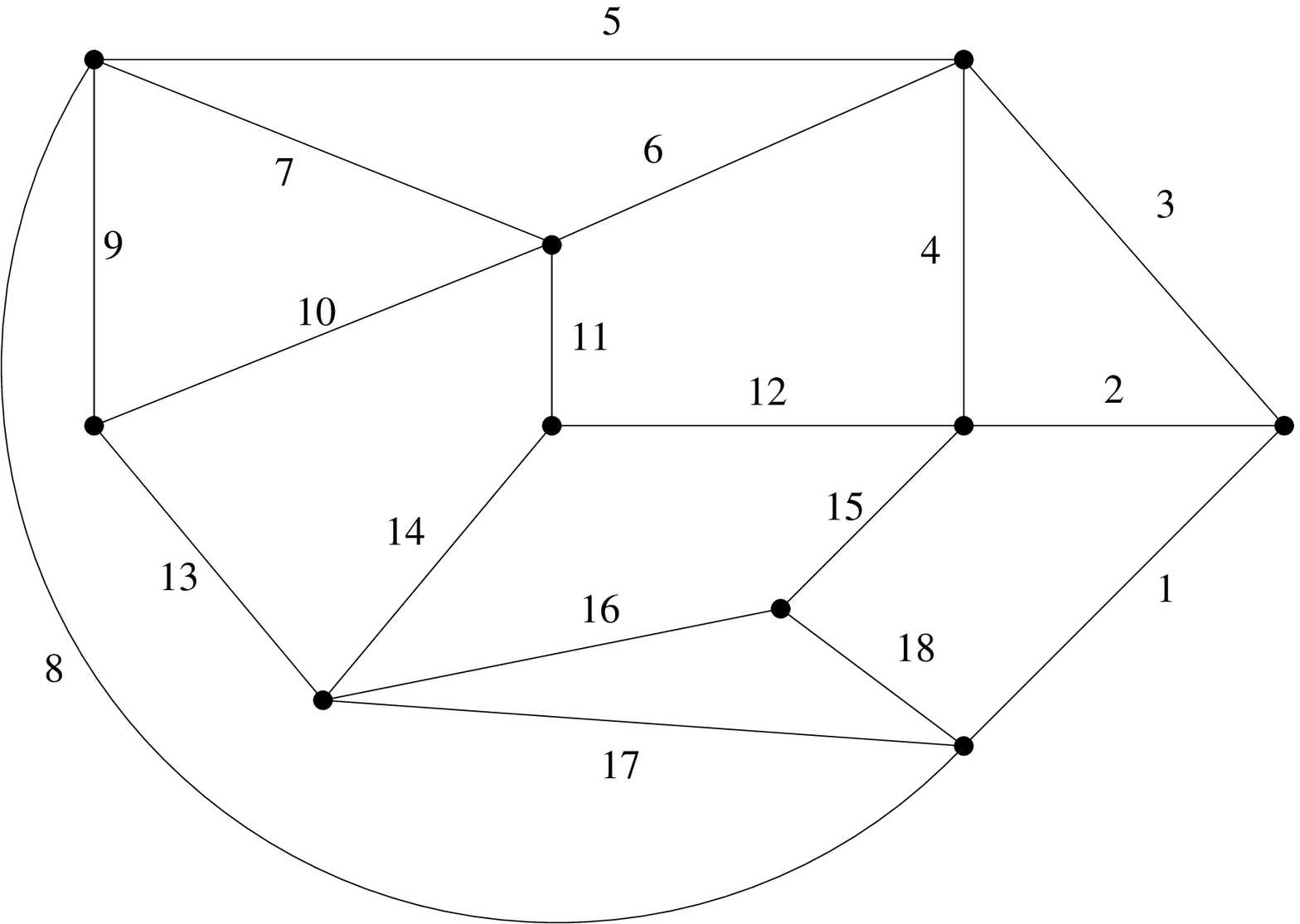}
      \put(-200,10){\large{$G_9$}}  
      \put(-140,88){{$t_1$}}  
      \put(-105,102){{$t_2$}}  
 \end{center}
\caption{A planar counter-example to Kontsevich's conjecture, with vertex width 4 (for the edge-ordering shown).}
\end{figure}
It contains a double triangle ($t_1$ and $t_2$), bounded by edges $5,6,7,9,10$. 
By applying a double-triangle reduction, the $c_2$-invariant
of this graph is equal to the $c_2$-invariant of a non-planar graph $G_9'$  at 8 loops. One verifies that the completion class of $G_9'$ is the same as that of $G_8$.
Thus, accepting the completion conjecture, we have
$c_2(G_9)_q \equiv 2-  [J]_q \mod q$ also.  In any case, a computer reduction of $G_9'$   (yielding a different quartic from $J$) confirms this prediction.

\section{A singular  K3 surface} \label{sect7}
 Consider the homogeneous polynomial of degree four
\begin{equation}\label{Fpolydef} F=b(a+c)(ac+bd) -ad(b+c)(c+d)
\end{equation}
which satisfies $F|_{d=1} = J$. One easily checks that it has six singular points
$$ e_1 = (0 : 0: 0: 1) \quad e_2 = (0 : 0 : 1 : 0) \quad e_3 = (0 : 1 : 0 : 0)$$
$$ e_4 = (1 : 0: 0: 0) \quad e_5 = (0 : 0 :-1 : 1) \quad e_6 = (1 : 1 :-1 : 1) $$
which are all of du Val type. Its minimal desingularization is obtained by blowing up the six points $e_1,\ldots, e_6$ and is therefore a K3 surface $X$. 
Since the Hodge numbers of a K3 satisfy $h^{1,1}=20$, and $h^{0,2}=h^{2,0}=1$, both $X$ and $V(F)\subset \Pro^3$ are not of Tate type and we can already conclude  by $(\ref{c2ofG8})$ that
the graph $G_8$  is   a counterexample to Kontsevich's conjecture by $(\ref{c2ofG8})$.

\subsection{The Picard lattice}
 We  determine the Picard lattice of $X$  as follows.
 It follows by inspection of $F$ that the following
lines  lie on $X$.
\begin{eqnarray}
\ell_1 : \,\, c=d=0 \qquad && \ell_{8} :\,\,  c= b+d=0\\ \nonumber 
\ell_2 :\,\,  b=d=0 \qquad&& \ell_{9} : \,\, b= c+d= 0\\ \nonumber 
\ell_3 : \,\, a=d=0 \qquad&& \ell_{10} : \,\, a-b=c+d=0 \\ \nonumber 
\ell_{4} : \,\, b=c=0\qquad && \ell_{11} : \,\, a=b=d \\ \nonumber 
\ell_{5} : \,\, a=c=0 \qquad&& \ell_{12}: \,\, a = b=-c \\ \nonumber 
\ell_{6} : \,\, a=b=0\qquad && \ell_{13} : \,\,  a=-c=d \\ \nonumber 
\ell_{7} : \,\, a+c=d =0 \,&& \ell_{14}: \,\,  a-d=b+c=0  \nonumber 
\end{eqnarray}  
Let $\ell_{15},\ldots, \ell_{20}$ denote the six exceptional divisors lying above the points $e_1,\ldots, e_6$.
Since these rational curves have self-intersection $-2$, one easily deduces the following intersection matrix, where the rows and columns correspond to $\ell_1,\ldots, \ell_{20}$.
$$\tiny 
 \left( \begin {array}{cccccccccccccccccccc} -\!2\!\!&0&0&0&0&0&0&0&0&1&0&0&0
&0&0&1&1&0&0&0\\\noalign{\medskip}0&-\!2\!\!&0&0&0&0&1&0&0&0&0&0&0&0&0&0&1&1
&0&0\\\noalign{\medskip}0&0&-\!2\!\!&0&0&0&0&0&0&0&0&0&0&1&0&1&0&1&0&0
\\\noalign{\medskip}0&0&0&-\!2\!\!&0&0&0&0&0&0&0&0&0&1&1&0&1&0&0&0
\\\noalign{\medskip}0&0&0&0&-\!2\!\!&0&0&1&0&0&0&0&0&0&1&1&0&0&0&0
\\\noalign{\medskip}0&0&0&0&0&-\!2\!\!&0&0&0&0&0&0&0&0&1&0&0&1&1&0
\\\noalign{\medskip}0&1&0&0&0&0&-\!2\!\!&0&0&0&0&1&0&0&0&1&0&0&0&0
\\\noalign{\medskip}0&0&0&0&1&0&0&-\!2\!\!&0&0&1&0&0&0&0&0&1&0&0&0
\\\noalign{\medskip}0&0&0&0&0&0&0&0&-\!2\!\!&0&0&0&1&0&0&0&1&0&1&0
\\\noalign{\medskip}1&0&0&0&0&0&0&0&0&-\!2\!\!&0&0&0&0&0&0&0&0&1&1
\\\noalign{\medskip}0&0&0&0&0&0&0&1&0&0&-\!2\!\!&0&0&0&0&0&0&1&0&1
\\\noalign{\medskip}0&0&0&0&0&0&1&0&0&0&0&-\!2\!\!&0&0&1&0&0&0&0&1
\\\noalign{\medskip}0&0&0&0&0&0&0&0&1&0&0&0&-\!2\!\!&0&0&1&0&0&0&1
\\\noalign{\medskip}0&0&1&1&0&0&0&0&0&0&0&0&0&-\!2\!\!&0&0&0&0&0&1
\\\noalign{\medskip}0&0&0&1&1&1&0&0&0&0&0&1&0&0&-\!2\!\!&0&0&0&0&0
\\\noalign{\medskip}1&0&1&0&1&0&1&0&0&0&0&0&1&0&0&-\!2\!\!&0&0&0&0
\\\noalign{\medskip}1&1&0&1&0&0&0&1&1&0&0&0&0&0&0&0&-\!2\!\!&0&0&0
\\\noalign{\medskip}0&1&1&0&0&1&0&0&0&0&1&0&0&0&0&0&0&-\!2\!\!&0&0
\\\noalign{\medskip}0&0&0&0&0&1&0&0&1&1&0&0&0&0&0&0&0&0&-\!2\!\!&0
\\\noalign{\medskip}0&0&0&0&0&0&0&0&0&1&1&1&1&1&0&0&0&0&0&-\!2\!\!
\end {array} \right)
$$ 
It has determinant $-7$.

  Since $7$ is prime, the lines $\ell_1,\ldots, \ell_{20}$ span the full 
N\'eron-Severi group. In particular, the rank of $X$ is 20 and so it defines a singular K3 surface. Since  $\Q(\sqrt{-7})$ has class number 1, 
$X$ corresponds to the unique singular K3  in the Shioda-Inose classification  \cite{SI} with
discriminant -7.
Now consider the elliptic curve $E=E_{49A1}$ with complex multiplication by $\Q(\sqrt{-7})$ which is given  by the affine model:
$$y^2+xy = x^3-x^2-2x-1\ . $$
The results of \cite{SI} imply that the graph of the complex multiplication in $E\times E$ gives rise to a decomposition of $\mathrm{Sym}^2 H^1(E)$ into two pieces, one of which is
$H^2_{tr}(X)$.
The results of Livn\'e \cite{Liv} allow one to conclude that the weight 3 modular form corresponding to $H^2_{tr}(X)$  is given by the symmetric square of the modular form of $E$.  
It is given explicitly by the following cusp form of weight 3 and level 7:
\begin{equation} \label{etas} \big( \eta(z)\eta(z^7)\big)^3
\end{equation}
where $\eta$ denotes the Dedekind eta function (first entry  of Table 2 in \cite{Schutt}).

\begin{remark}  Consider  Ramanujan's double theta function:
$$\theta(r,s)=\sum_{n=-\infty}^{\infty} r^{n(n+1)/2}s^{n(n-1)/2}$$
  and write $\theta_{a,b}(q) = \theta(-q^a,-q^b)$. Then, following \cite{Parry},   set
\begin{eqnarray} 
f_{49}(q)  &= &  \theta_{7,14}(q)^3 \big[ q \, \theta_{21,28}(q)  +  q^2 \theta_{14,35}(q) -q^4 \theta_{7,42}(q) \big] \nonumber \\
&= &q+q^2-q^4-3q^8-3q^9+4q^{11}-q^{16}-3q^{18}+4q^{22}+8q^{23}+\ldots   \nonumber 
\end{eqnarray}
which spans the one-dimensional space of  newforms of level 49 and weight 2 (see  also \cite{Sloane}).
   If $a_{p^n}$ denotes the coefficient of $q^{p^n}$ in $f_{49}(q)$, 
one knows that  the number of points of $E$ over $\F_{p^n}$ is  $p^n+1-a_{p^n}$. One can show that:
\begin{equation}\label{apcases}
 a_p = \left\{  
\begin{array}{ll}
  0   &  \hbox{ if }\quad  p \equiv 0,3,5, 6 \mod 7 \ ,  \\
   \pm a  & \hbox{ where }  4p=a^2+7 b^2 \hbox{ where }    a,b\in \Z, \hbox{ if } p \equiv 1,2,4 \mod 7  \ . 
\end{array}   \right.
 \end{equation} 
\end{remark}

Let  $b_{p^n}$ denote the coefficient of $z^{p^n}$ in $(\ref{etas})$.  Modulo $p$, we simply have 
\begin{equation}\label{a2b}
a^2_{p^n} \equiv b_{p^n} \mod p\ .
\end{equation}

\begin{thm} Let $G_8$ be the 8-loop non-planar graph defined in \S\ref{sect8loopcounterNP}, and figure 8. Then the number of points of the affine graph hypersurface $X_{G_8}$
over $\F_{p^n}$ satisfies:
\begin{equation}\label{finalXG8expression} 
 [X_{G_8}]_{{p^n}} \equiv  - a^2_{p^n} p^{2n}  \equiv   -b_{p^n} p^{2n} \pmod  {p^{2n+1}} \ . 
 \end{equation} 
\end{thm}
\begin{proof} 
Let $q=p^n$. We have  $[X_{G_8}] \equiv c_2(G_8)\, q^2 \mod q^3$. Equation $(\ref{c2ofG8})$ states that  $c_2(G_8)\equiv 2-[J]_q \mod q$.  Passing to the 
homogeneous version $(\ref{Fpolydef})$, one verifies that $[J]_q\equiv 2 - [F]_q \mod q$.
Finally, the above discussion  and equation $(\ref{a2b})$ shows that  $[F]_q \equiv - a^2_q \equiv - b_q \mod p$.  Therefore $[X_{G_8}]$  modulo  $q^2 p$ is  given by 
$(\ref{finalXG8expression})$.
\end{proof}
Consider the product of all finite fields $\F_p$ where $p$ is prime:
$$A = \F_2\times \F_3 \times \F_5 \times  \ldots $$
and define the total $c_2$-invariant of 	a graph $G$ to be 
$$\widetilde{c}_2(G) =   (c_2(G)_2,c_2(G)_3, c_2(G)_5,\ldots   ) \in A \  ,$$
where we identify $\Z/p\Z$ with $\F_p$.  Let $\pi: \Z \rightarrow A$ denote the map whose $p^\mathrm{th}$ component is $ n \mapsto n \mod p$. 
It follows from $(\ref{apcases})$ that 
$$\widetilde{c}_2(G_8) \notin \pi(K) $$
for all bounded sets $K\subset \Z$, since in the opposite case, all primes $p$ congruent to $1,2,4$ mod $7$ would satisfy
$4p \in - K+  \{7 \,b^2: b\in \Z\}$. Since $K$ is finite, this would  contradict the prime number theorem. Therefore $c_2(G_8)$ is not (quasi-)constant, and therefore the graph $G_8$ is  a counter-example to Kontsevich's conjecture in the strongest possible sense.

Assuming the  completion conjecture \ref{completionconjecture}, or by the double-triangle theorem and the computer calculation in \cite{SchnetzFq},  the  graph $G_9$ of
\S\ref{sectplanarcounter} has  exactly the same property, and yields a planar counter-example at 9 loops.

\subsection{Discussion}
The prevalence  of multiple zeta values in Feynman integral computations  at low loop orders led Kontsevich to conjecture that the Euler characteristics
of graph hypersurfaces were of mixed Tate type. This was shown to be generically false by Belkale and Brosnan, but despite this 
cautionary result,  the following questions about the arithmetic nature of $\phi^4$ theory remained open:
\begin{enumerate}
\item Even though general graphs have non-Tate Euler characteristics, it could be that graphs coming from physically relevant theories
are still of Tate type (the counter-examples of \cite{BB} have unphysical 
vertex degrees).
\item It could be that the counter-examples  occur at very high loop order rendering them physically less relevant.
\item Failing (1) and (2), it could still be the case that planar graphs have Tate Euler characteristics,  i.e., all non-Tate counter-examples can  be characterized by having a
high genus or crossing number.
\item Even though the Euler characteristics are non-Tate, it could be that the piece of the graph motive which carries the period is always mixed Tate.
\end{enumerate}
Our counter-examples   show that (1), (2) and (3) are false.  Point (4) is more subtle. 
However, it follows from the original interpretation of the denominator reduction in \cite{BrFeyn} that the $c_2$-invariant of a graph should correspond
to the `framing' on $M_G$, i.e., the smallest subquotient motive of $M_G$ which is spanned by the integrand of $(\ref{IG})$.  
This makes it very probable  that (4) is  false too. In this case, one is led to expect the residue $I_G$ (equation $(\ref{IG})$) to be transcendental over the ring generated by
multiple zeta values over $\Q$. Indeed, a likely candidate for the periods of the  counter-examples $G_8, G_9$ might   come from the periods of the  motivic fundamental 
group of the   elliptic curve $E_{49A1}$ with punctures.

Finally, it should be emphasized that the residues $I_G$ of primitive graphs in $\phi^4$ are renormalization-scheme independent, and universal in the sense that any
quantum field theory in 4 space-time dimensions  will only affect the numerator, and not the denominator, of the corresponding parametric integral representation (barring infra-red divergences). Since the motive $M_G$ only
depends on the denominators, one can reasonably expect that such non-mixed Tate  phenomena will propagate into most  renormalizable  massless  quantum field theories
with a four-valent vertex at sufficiently high loop orders.
\vskip2ex

\bibliographystyle{plain}
\bibliography{main}

\begin{thebibliography}{99}

\bibitem{AM} {\bf P. Aluffi, M. Marcolli}: {\it Feynman motives and deletion-contraction relations},
In: Topology of Algebraic Varieties and Singularities,  Contemporary Mathematics 538, AMS (2011).

\bibitem{BB} {\bf P. Belkale, P. Brosnan}: {\it Matroids, motives, and a conjecture of Kontsevich}, Duke Math.\ J. 116, no.\ 1, 147-188 (2003).

\bibitem{BK}{\bf D. Broadhurst, D.  Kreimer}: {\it  Knots and numbers in $\phi^4$ theory to 7 loops and beyond}, Int.\ J. Mod.\ Phys.\ C 6, 519 (1995).

\bibitem{BlochJ} {\bf S. Bloch}:  {\it Motives associated to sums of  graphs},  \url{http://www.math.uchicago.edu/~bloch/publications.html}.

\bibitem{B-E-K} {\bf S. Bloch, H. Esnault, D. Kreimer}:  {\it On motives associated to graph polynomials},
Comm.\ Math.\ Phys.\ 267, no.\ 1, 181-225 (2006).

\bibitem{BrAn} {\bf F. Brown}:  {\it Mixed Tate motives over $\Z$},  to appear in Annals of Math., (2011).


\bibitem{BrCMP} {\bf F. Brown}:  {\it The massless higher-loop two-point function}, Comm.\ in Math.\ Physics 287, no.\ 3, 925-958 (2009).

\bibitem{BrFeyn}{\bf F.  Brown}:  {\it On the periods of some Feynman graphs}, \url{arXiv:0910.0114v1} (2009).

\bibitem{WD} {\bf  F. Brown, K. Yeats}: {\it Spanning forest polynomials and the transcendental weight of Feynman graphs}, 
Comm. in Math. Physics,  vol. 301,  pp. 357-382 (2011)
 
\bibitem{CY} {\bf F. Chung, C. Yang}: {\it On Polynomials of Spanning trees}, Annals of Combinatorics 4, 13-25 (2000).

\bibitem{Doryn} {\bf D. Doryn}: {\it Cohomology of graph hypersurfaces associated to certain Feynman graphs}, 
 \url{arXiv:0811.0402} (2008).

\bibitem{DorynPoints} {\bf D. Doryn}: {\it On one example and one counterexample in counting rational points on graph hypersurfaces}, 
 \url{arXiv:1006.3533v1} (2010).

\bibitem{EV} {\bf H. Esnault, E. Viehweg} : {\it On a rationality question in the Grothendieck ring of varieties}, \url{arXiv:0908.2251v3} (2009).

\bibitem{Kir}{\bf G. Kirchhoff}: {\it Ueber die Aufl\"osung der Gleichungen, auf welche man bei der Untersuchung der linearen Vertheilung galvanischer Str\"ome geh\"uhrt wird},
 Annalen der Physik und Chemie 72, no.\ 12, 497-508 (1847).

\bibitem{Liv}{\bf R. Livn\'e}: {\it Motivic Orthogonal Two-Dimensional Representations of $\hbox{Gal}(\overline{Q}/Q)$}, Israel  journal of mathematics 92, 149-156 (1995).

\bibitem{Parry}{\bf W. Parry}:  {\it  A negative result on the representation of modular forms by theta series}, J.\ Reine Angew.\ Math.\ 310, 151-170  (1979).

\bibitem{SchnetzCensus} {\bf O. Schnetz}: {\it Quantum periods: A census of $\phi^4$ transcendentals}, Jour.\ Numb.\ Theory
and Phys.\ 4 no.\ 1, 1-48 (2010).

\bibitem{SchnetzFq} {\bf O. Schnetz}: {\it Quantum field theory over $F_q$}, The Electronic Jour.\ of Comb.\ 18, \#P102 (2011).

\bibitem{Schutt} {\bf M. Sch\"utt}: {\it CM newforms with rational coefficients}, Ramanujan Journal 19, 187-205 (2009).

\bibitem{SI} {\bf T. Shioda, H.  Inose}: {\it On Singular K3 Surfaces},  in  Complex Analysis
and Algebraic Geometry, ed. W. Baily and T. Shioda,  Iwanami Shoten, Tokyo, pp. 119-136 (1977).

\bibitem{Sloane} {\bf N. Sloane}: {\it Online encyclopedia of integer sequences}, A140686. 

\bibitem{Se} {\bf J.-P.  Serre}: {\it Cours d'arithm\'etique}, Presse Universitaire de France, (1977).

\bibitem{Sta} {\bf  R. P. Stanley}: {\it Spanning trees and a conjecture of
Kontsevich}, Ann.\ Comb.\ 2, no.\ 4, 351-363 (1998).

\bibitem{Stem} {\bf J. Stembridge}: {\it Counting points on varieties over finite fields related to a
conjecture of Kontsevich},  Ann.\ Combin.\ 2,  365--385 (1998).

\bibitem{Wein} {\bf S. Weinberg}: {\it High-Energy Behavior in Quantum Field Theory},  Phys.\ Rev.\ 118, no.\ 3 838--849 (1960).

\end{thebibliography}

\end{document}